\documentclass[12pt, a4paper]{article}

\usepackage[vmargin = 30mm, hmargin = 25mm, a4paper]{geometry}
\usepackage[utf8]{inputenc}
\usepackage{lmodern}
\usepackage{amsmath,amsfonts,amsthm}
\usepackage{authblk}
\usepackage{caption}
\usepackage{float}
\usepackage{graphicx}
\usepackage[numbers]{natbib}
\usepackage{parskip}
\usepackage[section]{placeins}
\usepackage{setspace}
\usepackage{subcaption}
\usepackage[dvipsnames]{xcolor}
\usepackage{hyperref}
\usepackage[noabbrev]{cleveref}

\usepackage[mathlines]{lineno}
\newcommand*\linenomathpatch[1]{%
  \cspreto{#1}{\linenomath}%
  \cspreto{#1*}{\linenomath}%
  \csappto{end#1}{\endlinenomath}%
  \csappto{end#1*}{\endlinenomath}%
}
\newcommand*\linenomathpatchAMS[1]{%
  \cspreto{#1}{\linenomathAMS}%
  \cspreto{#1*}{\linenomathAMS}%
  \csappto{end#1}{\endlinenomath}%
  \csappto{end#1*}{\endlinenomath}%
}

\expandafter\ifx\linenomath\linenomathWithnumbers
  \let\linenomathAMS\linenomathWithnumbers
  \patchcmd\linenomathAMS{\advance\postdisplaypenalty\linenopenalty}{}{}{}
\else
  \let\linenomathAMS\linenomathNonumbers
\fi

\linenomathpatch{equation}
\linenomathpatchAMS{gather}
\linenomathpatchAMS{multline}
\linenomathpatchAMS{align}
\linenomathpatchAMS{alignat}
\linenomathpatchAMS{flalign}

\makeatletter
\patchcmd{\mmeasure@}{\measuring@true}{
  \measuring@true
  \ifnum-\linenopenaltypar>\interdisplaylinepenalty
    \advance\interdisplaylinepenalty-\linenopenalty
  \fi
  }{}{}
\makeatother

\crefname{equation}{}{}

\newcommand{\e}{{\rm e}}
\newcommand{\ShockMin}{u_l}
\newcommand{\ShockMax}{u_r}

\renewcommand\Re{\operatorname{Re}}

\newcommand{\ve}{\varepsilon}
\newcommand{\cL}{\mathcal{L}}
\newcommand{\cM}{\mathcal{M}}
\newcommand{\C}{\mathbb{C}}
\newcommand{\R}{\mathbb{R}}

\DeclareRobustCommand{\d}{\relax\ifmmode\mathrm{d}\else\expandafter\@@d\fi}
\DeclareRobustCommand{\D}{\relax\ifmmode\mathrm{D}\else\expandafter\@@d\fi}
\DeclareRobustCommand{\e}{\relax\ifmmode\mathrm{e}\else\error\fi}
\DeclareRobustCommand{\i}{\relax\ifmmode\mathrm{i}\else\expandafter\@@i\fi}
\DeclareRobustCommand{\uppartial}{\text{\rotatebox[origin=t]{20}{\scalebox{0.95}[1]{\(\partial\)}}}\hspace{-1pt}}

\newcommand{\fd}[2]{\frac{\d #1}{\d #2}}

\newcommand{\fdn}[3]{\frac{\d^#3 #1}{\d #2^#3}}
\newcommand{\pd}[2]{\frac{\uppartial #1}{\uppartial #2}}

\newcommand{\pdn}[3]{\frac{\uppartial^#3 #1}{\uppartial #2^#3}}

\newcommand{\df}[1]{\,\d#1}

\DeclareSymbolFont{vectors}{OML}{lmm}{b}{it}
\DeclareSymbolFontAlphabet{\mathvec}{vectors} % Vectors

\newtheorem{theorem}{Theorem}[section]

\theoremstyle{definition}
\newtheorem*{remark}{Remark}

\def\arev{\color{black}} % Denote revision
\def\brev{\color{black}} % Denote revision
\def\endrev{\color{black}}

\title{Analytic shock-fronted solutions to a reaction--diffusion equation with negative diffusivity}
\author[1]{Thomas Miller}
\author[1]{Alexander K. Y. Tam}
\author[2]{Robert Marangell}
\author[2]{Martin Wechselberger}
\author[1,*]{Bronwyn H. Bradshaw-Hajek}
\affil[1]{UniSA STEM, The University of South Australia, Mawson Lakes SA 5095, Australia}
\affil[2]{School of Mathematics and Statistics, The University of Sydney, Sydney NSW 2006, Australia}
\affil[*]{Corresponding author: Bronwyn H. Bradshaw-Hajek, \href{mailto:bronwyn.hajek@unisa.edu.au}{bronwyn.hajek@unisa.edu.au}}

\begin{document}
\maketitle     
\doublespacing
% \linenumbers
%
%%%%%%%%%%%%%%%%%%%%
%%%%% Abstract %%%%%
%%%%%%%%%%%%%%%%%%%%
%
\section*{Abstract}
Reaction--diffusion equations (RDEs) model the spatiotemporal evolution of a density field \(u(\mathvec{x},t)\) according to diffusion and net local changes. Usually, the diffusivity is positive for all values of \(u,\) which causes the density to disperse. However, RDEs with \brev partially \endrev negative diffusivity can model aggregation, which is the preferred behaviour in some circumstances. In this paper, we consider a nonlinear RDE with quadratic diffusivity \(D(u) = (u - a)(u - b)\) that is negative for \(u\in(a,b)\). We use a nonclassical symmetry to construct analytic receding time-dependent, colliding wave, and receding travelling wave solutions. \brev These solutions are multi-valued, \endrev and we convert them to single-valued solutions by inserting a shock. We examine properties of these analytic solutions including their Stefan-like boundary condition, and perform a phase plane analysis. We also investigate the spectral stability of the \(u = 0\) and \(u = 1\) constant solutions, and prove for certain \(a\) and \(b\) that receding travelling waves are spectrally stable. Additionally, we introduce \arev a \endrev new shock condition where the diffusivity and flux are continuous across the shock. For diffusivity symmetric about the midpoint of its zeros, this condition recovers the well-known equal-area rule, but for non-symmetric diffusivity it results in a different shock position.

\paragraph{Keywords:} nonclassical symmetry, spectral stability, phase plane analysis, aggregation, travelling wave, equal-area rule
%
%%%%%%%%%%%%%%%%%%%%%%%%%%%
%%%%% 1. Introduction %%%%%
%%%%%%%%%%%%%%%%%%%%%%%%%%%
%
\section{Introduction}\label{section:Introduction}%
Reaction--diffusion equations (RDEs) are partial differential equation (PDE) models that have countless applications, including in wound healing~\cite{Models_of_epidermal:sherratt1990}, biological invasion~\cite{Allee_dynamics:lewis1993}, and movement through porous media~\cite{Symmetry_solutions:moitsheki2005}. The simplest reaction--diffusion models consider the evolution of one population in one spatial dimension. Mathematically, these models have the general form
\begin{equation}
    \label{equation:Reaction_diffusion}%
    \pd{u}{t} = \pd{}{x}\left(D(u)\pd{u}{x}\right) + R(u), 
\end{equation}    
where \(u(x,t)\) is the density as a function of space, \(x,\) and time, \(t.\) In~\Cref{equation:Reaction_diffusion}, the symbol \(D(u)\) is the \brev diffusivity \endrev and the function \(R(u)\) is the reaction term capturing local changes in density.  In standard linear Fickian diffusion, \(D(u) = D,\) where \(D > 0\) is a constant. However, many \brev systems \endrev do not undergo linear diffusion, and in these scenarios, \(D(u)\) is non-constant. \brev In general, we refer to non-constant \(D(u)\) \endrev as nonlinear diffusion.
        
Reaction--diffusion models with nonlinear diffusion commonly have \(D(u) \geq 0\) for all feasible values of \(u.\) Positive diffusivity results in the population dispersing from regions of high density to regions of low density. However, many natural \brev processes \endrev involve aggregation, whereby members of the population move up density gradients from regions of low density to regions of high density. For example, a near-universal phenomenon in ecology is animal populations aggregating for social reasons including mating, prey detection, or predator avoidance~\citep{Modelling_social:grunbaum1994}. \arev Typical models for aggregation involve members of the population responding to chemical stimuli~\citep{Keller1970}. However, another \endrev way to model aggregation is to allow the diffusivity \(D(u)\) to become negative for some values of \(u\). Negative diffusivity can arise in the continuum limit of discrete models with social aggregation~\citep{Model_for_spreading:popescu2004, Population_consequences:turchin1989}. In~\citep{Model_for_spreading:popescu2004}, the strength of inter-particle attraction influences \brev the regions in which the diffusion is \endrev negative. In~\citep{Population_consequences:turchin1989}, by considering a random walk with attraction among individuals \citet{Population_consequences:turchin1989} derived a model with a diffusivity that could become negative. \citet{Co-operation_competition:johnston2017} showed that an agent-based model containing both isolated and group-associated individuals which can have different rates of birth, death and movement gives rise to an RDE model with a region of negative diffusivity. In this work, we focus on models with \(D(a) = D(b) = 0,\) with \(D(u) < 0\) for \(a < u < b,\) where \(a\) and \(b\) are constants such that \(0 < a < b < 1,\) with \(D(u) \geq 0\) otherwise. \arev Although it remains to be demonstrated that negative diffusivity can describe aggregation in a real system, negative diffusivity provides a way to generate aggregation using a single PDE, without requiring chemotaxis. \endrev
        
\arev Allowing negative diffusivity \(D(u)\) causes ill-posedness in Cauchy initial-value problems involving the nonlinear diffusion equation \(u_t = \left(D(u) u_x\right)_x\)~\citep{Models_for_mutual:alt1985, The_role_of_diffusion:aronson1985}. It also causes solutions for \(u\) to have very steep gradients~\cite{Population_consequences:turchin1989}\endrev, and solutions can even become multi-valued. Inserting a jump discontinuity (or shock) in a multi-valued solution enables a single-valued solution to be obtained. A shock involves an instantaneous (with respect to \(x\)) jump from one density, \(\ShockMax\), to another, \(\ShockMin\), such that \(u\) never takes the values between \(\ShockMax \geq b\) and \(\ShockMin \leq a\). Consequently, both the region of negative diffusivity \(u \in (a,b)\), and the region where the solution is multi-valued are avoided altogether. Inserting a shock is an established technique for obtaining single-valued solutions to hyperbolic equations of the form \(\rho_t + c(\rho)\rho_x = 0\)~\cite{Linear_and_nonlinear_waves:whitham2011}. \citet{Shocks_in_nonlinear_diffusion:witelski1995} expanded this technique to construct a solution for a diffusion equation with a region of negative diffusivity. Furthermore, some reaction--diffusion models with a region of negative diffusivity admit travelling wave solutions where the density profile moves at constant speed while maintaining its shape~\cite{Travelling_waves_in_some:ferracuti2009, Front_propagation:kuzmin2011, Travelling_wave_solutions:li2020, Shock-fronted:li2021, Diffusion_aggregation:maini2006}. Of particular interest to this work, \citet{Shock-fronted:li2021} proved the existence of shock-fronted travelling waves for one such reaction--diffusion model.
        
We use a nonclassical Lie symmetry to construct a multi-valued analytic solution to the RDE~\Cref{equation:Reaction_diffusion} with a region of negative diffusivity. Lie symmetry analysis was first proposed by Sophus Lie in the late 19th century~\cite{Lie_article} and involves finding invariant quantities in an equation and using them to construct solutions (see~\cite{Elementary_Lie:ibragimov1999, Applications_of_Lie:olver1993} for more detail). The nonclassical (or Q-conditional) method is an extension to Lie symmetry analysis introduced by~\citet{The_general_similarity:bluman1969}, where in addition to requiring invariance of the PDE the invariant surface condition must also be satisfied. Nonclassical symmetry analysis can sometimes uncover symmetries not found by the classical method. \brev While boundary and initial conditions are important for determining the behaviour of a system, they are often not included in the initial analysis because calculation of the symmetries becomes overly restrictive~\cite{Goard2008}. We follow this strategy here and apply the symmetry method to the governing equation only, revisiting the initial and boundary conditions after a solution has been found. \endrev
% Mention that the nonclassical symmetry analysis usually does not involve specifying any IC or BC at the start - we find the possible analytic solutions then choose those that have sensible properties.
            
Here we use a special nonclassical symmetry found independently by~\citet{Nonclassical_symmetries:arrigo1995} and~\citet{Nonclassical_symmetry:goard1996} to construct analytic solutions for an arbitrary diffusivity provided that the reaction and diffusion terms satisfy the relationship 
\begin{equation}
    R(u) = \left(\frac{A}{D(u)} + \kappa\right)\int_{u^*}^u D(u') \df{u'},
    \label{equation:Reaction}
\end{equation} 
where \(u^*\) is a zero of the reaction term, and \(A\) and \(\kappa\) are free  parameters that arise from the symmetry analysis. Generally this reaction term contains singularities \brev (simple poles) \endrev when \(D(u) = 0\) (that is, at \(u = a\) and \(u = b\)), however inserting a shock avoids these singularities. \arev For the applications considered here, the case that \(u^*\) is a zero of \(D(u)\) as well as a zero of \(R(u),\) is not considered. \endrev This symmetry \brev transforms \endrev the RDE~\Cref{equation:Reaction_diffusion} to the \brev one-dimensional \endrev spatial Helmholtz equation~\cite{Conditionally_Integrable:broadbridge2023},
\begin{equation}
    \label{equation:Helmholtz}%
    \fdn{\Psi}{x}{2} + \kappa \Psi = 0.
\end{equation}
The new dependent variable \(\Psi(x)\) is related to the density \(u(x,t)\) via the Kirchhoff transformation
\begin{equation}
    \label{equation:Transformation}%
    \Phi(u) = \int_{u^*}^u D(\bar{u}) \df{\bar{u}} = e^{At}\Psi(x),
\end{equation}
where \(\Phi(u)\) is the Kirchhoff variable, sometimes known as the diffusive flux or potential. After solving the Helmholtz equation~\Cref{equation:Helmholtz}, inverting the Kirchhoff transformation~\Cref{equation:Transformation} recovers the solution for \(u(x,t)\). We restrict our attention to the choice \(\kappa = -k^2 < 0\), \brev so that \endrev the solution to the Helmholtz equation~\Cref{equation:Helmholtz} is in terms of exponential functions,
\begin{equation}
    \Psi(x) = c_1 e^{kx} + c_2 e^{-kx}.
    \label{equation:Psi}        
\end{equation}
In summary, we will specify a diffusivity with a region of negative diffusion and calculate the related reaction term~\Cref{equation:Reaction}. We then construct analytic solutions by inverting the transformation~\Cref{equation:Transformation} using the exact form of $\Psi(x)$~\eqref{equation:Psi}.        
                
In~\Cref{section:ExampleSolution}, we investigate a model where \(D(u)\) is a quadratic function with two zeros on \(u \in (0,1),\) such that there is a region \(u \in (a,b),\) where \(0 < a < b < 1,\) \brev with \endrev \(D(u) < 0.\) We consider three types of analytic solutions. The first two are time-dependent solutions, which can represent both a receding front with non-constant speed and a pair of colliding waves. The third is a receding sharp-fronted travelling wave, where the solution moves leftward at constant speed. \brev All three solution types are sharp-fronted and satisfy a Stefan-like moving boundary condition due to the requirement that the density remains non-negative \cite{Compactly_Supported_Solutions:edwards2018, Revisiting_the_Fisher:el-hachem2019}. \endrev Since the solutions with negative diffusivity are multi-valued, we create shock-fronted solutions by inserting a jump in \(u\) around values for which \(D(u) < 0\), \brev so that \endrev the shock \arev conserves \(\Phi'(u) = D(u)\) across the shock or \(\Phi'(u)\) is continuous rule. \endrev After developing the shock-fronted solutions, we analyse the stability of the constant and travelling wave solutions. \Cref{section:Nonsymmetric_diffusion} examines the relationship between the \arev\(\Phi'(u)\) is continuous rule \endrev and the more common equal area in \(\Phi(u)\) rule \cite{Front_migration:pego1989,Shocks_in_nonlinear_diffusion:witelski1995}. We provide discussion and concluding remarks in~\Cref{section:Discussion}. The existence of an analytic solution obtained using nonclassical symmetry analysis underpins the analysis in each section of the paper.
%
%%%%%%%%%%%%%%%%%%%%%%%%%%%%%%%%%%%%%%%%%%%%%%%%%%%%%%%%%%%
%%%%% 2: Analytic solutions for quadratic diffusivity %%%%%
%%%%%%%%%%%%%%%%%%%%%%%%%%%%%%%%%%%%%%%%%%%%%%%%%%%%%%%%%%%
%
\section{Analytic travelling wave and time-dependent solutions for quadratic diffusivity} \label{section:ExampleSolution}
We consider three analytic solutions for the RDE~\Cref{equation:Reaction_diffusion} with a quadratic diffusivity that is negative only for values of \(u\) between the two roots of the quadratic. This diffusivity takes the form
\begin{equation}        
    D(u) = \left(u - a\right)\left(u - b\right), \label{equation:Diffusion}
\end{equation}
where \(a,b \in (0,1)\) and we impose \(a < b\). \Citet{Co-operation_competition:johnston2017} derived a reaction--diffusion model with diffusivity~\cref{equation:Diffusion} in the continuum limit of a discrete model where isolated and group-associated individuals had different rates of birth, death and movement. \citet{Population_consequences:turchin1989} also derived a diffusion equation with a negative diffusivity from a discrete model where individuals move towards regions of higher population density~\citep{Population_consequences:turchin1989}. \Citet{Travelling_wave_solutions:li2020, Shock-fronted:li2021} then showed that travelling wave solutions exist for the RDE with quadratic diffusivity~\cref{equation:Diffusion}. We extend these works by using a nonclassical symmetry to construct solutions with shocks, and analyse the spectral stability of travelling wave solutions.                 
        
\begin{figure}[htbp!]
    \centering
    \subcaptionbox{\label{figures:RDP_a}}{\includegraphics[width=0.45\linewidth]{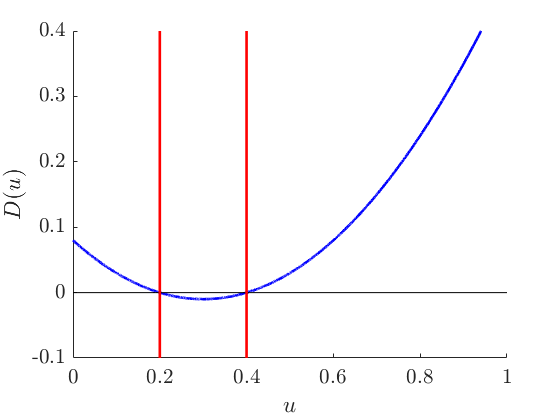}} 
    \subcaptionbox{\label{figures:RDP_b}}{\includegraphics[width=0.45\linewidth]{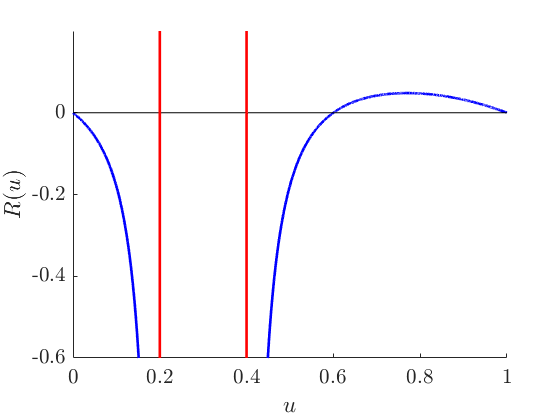}}
    \subcaptionbox{\label{figures:RDP_c}}{\includegraphics[width=0.45\linewidth]{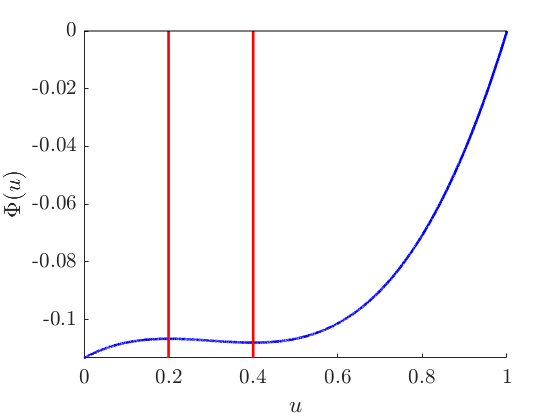}}
    \caption{Example nonlinear (a) diffusivity, (b) reaction function, and (c) flux potential. Here, \(\kappa = -1\), \(A = 0.08\), \(u^* = 1\), \(a = 0.2\) and \(b = 0.4\). The vertical lines (red) show the location of the zeros of the diffusivity. The reaction term has zeros at \(u = 0\), \(u = a + b = 0.6\), and \(u = 1.\) Between the vertical (red) lines \(R(u)\) becomes large and positive.}              
    \label{figures:RDP}
\end{figure}

\Cref{figures:RDP_a} shows an example quadratic diffusivity of the form~\cref{equation:Diffusion}, with \(a = 0.2\) and \(b = 0.4\). To use the special nonclassical symmetry, we obtain the corresponding reaction term using~\Cref{equation:Reaction}. By choosing \(u^* = 1,\) we fix one of the zeros of the reaction term to be at \(u = 1,\) that is, \(\Phi(1) = 0,\) and so \(R(1) = 0.\) Since \(A\) is a free parameter, we choose \(A = -\kappa D(0),\) to enforce \(R(0)=0.\) The reaction term \(R(u)\)~\Cref{equation:Reaction} has a third zero at \(u = a + b\). Therefore, there are two types of reaction terms possible for this quadratic diffusivity. If \(a + b \geq 1,\) the reaction term is typically negative for all values of \(u\) in the range of interest. Since a negative reaction term would cause a decrease in population for all densities, we do not investigate this case further. \brev If \(a + b < 1\), the reaction term is typically negative for lower densities and positive for higher densities. This scenario is loosely analogous to cubic reaction term commonly used to model changes in gene frequencies \cite{Huxley_and_Fisher:broadbridge2002}, or population dynamics due to the Allee effect \cite{Inverse_density:courchamp1999, Consequences_of_the_Allee:stephens1999}. \Cref{figures:RDP_b} shows the shape of this reaction term. \endrev The reaction term has singularities at \(a\) and \(b\), as the solid vertical lines in~\Cref{figures:RDP_b} indicate. Between the vertical lines the reaction term becomes large and positive. However, the shock solutions we construct avoid the \(u\) values in this region.
        
We construct an analytic solution to \Cref{equation:Reaction_diffusion} with diffusivity~\Cref{equation:Diffusion} and reaction term~\Cref{equation:Reaction} by solving the Helmholtz equation \Cref{equation:Helmholtz} for $\Psi(x)$, and inverting the Kirchhoff transformation \Cref{equation:Transformation} to find \(u(x,t)\). Since the diffusivity is quadratic, the Kirchhoff variable \(\Phi(u)\) will be cubic in \(u(x,t)\). The Kirchhoff transformation~\Cref{equation:Transformation} then determines the solution implicitly. That is,
\begin{equation}
    \Phi(u)=\frac{1}{6}\left(u - 1\right)\left(2u^2 + \left(2 - 3a - 3b\right)u + 2 - 3a - 3b + 6ab\right) = e^{At}\left(c_1\e^{kx} + c_2 \e^{-kx}\right),
    \label{equation:Implicit_solution}
\end{equation}
where \(c_1\) and \(c_2\) are constants of integration whose values will be governed by the boundary conditions, and where we have chosen \(\kappa = -k^2 < 0\). Some sections of the solution for \(u(x,t)\) are multi-valued, because \(\Phi(u)\) is a cubic function. Since \(A=-\kappa D(0)>0\), we require \(\Psi(x)<0\) so that the solution does not approach infinity as \(t\to\infty\). As a consequence, we require \(\Phi(u)<0\) for \(u\in[0,1)\), so that we must have \(b<(a+2)/3\). \Cref{figures:RDP_c} shows \(\Phi(u)\) which is increasing except where \(u\in(a,b)\) where the diffusivity is negative. Consequently the solution is multi-valued for values of \(u\) surrounding this region.   
%
%%%%%%%%%%%%%%%%%%%%%%%%%%%%%%%%%%%%%%%%%%%%%%%%
%%%%% Constructing single-valued solutions %%%%%
%%%%%%%%%%%%%%%%%%%%%%%%%%%%%%%%%%%%%%%%%%%%%%%%
%
\subsection{Constructing single-valued solutions by inserting shock discontinuities}\label{section:Shock}
\begin{figure}[htbp!]
    \centering
    \subcaptionbox{\label{figures:RDP_shock_a}}{\includegraphics[width=0.45\linewidth]{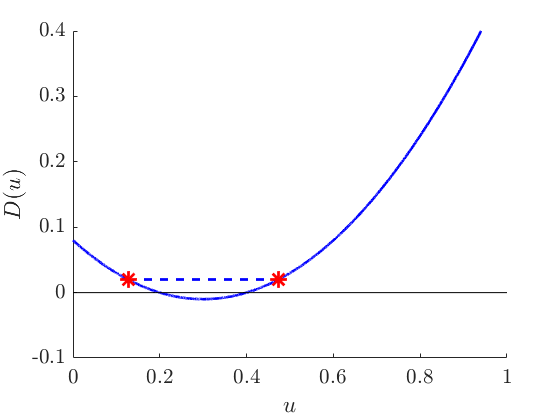}} 
    \subcaptionbox{\label{figures:RDP_shock_b}}{\includegraphics[width=0.45\linewidth]{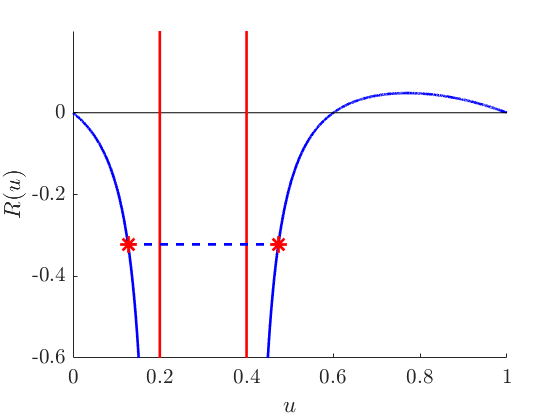}}            
    \subcaptionbox{\label{figures:RDP_shock_c}}{\includegraphics[width=0.45\linewidth]{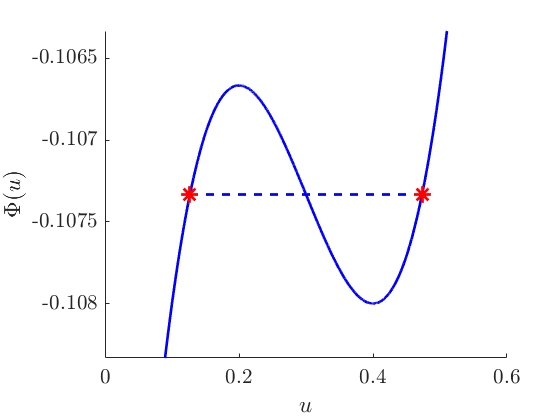}}
    \caption{Nonlinear (a) diffusivity, (b) reaction function, and (c) flux potential around the shock. Here, \(\kappa = -1\), \(A = 0.08\), \(u^* = 1\), \(a = 0.2\) and \(b = 0.4\). The vertical lines (red) show the location of the zeros of the diffusivity. Horizontal dashed lines indicate the transition across the shock, and the stars (red) indicate the end points of the shock. The reaction term is zero at \(u = 0\), \(u = a + b = 0.6\), \(u = 1\).}
    \label{figures:RDP_shock}
\end{figure}

With quadratic diffusivity \eqref{equation:Diffusion}, the solutions~\cref{equation:Implicit_solution} obtained using the nonclassical symmetry solutions described here are multi-valued with three branches in a region centred about the points where the diffusivity is zero, \emph{i.e.} where \(u = a\) and \(u = b.\) That is, the solution is multi-valued when \(u(x,t)\in\left[(3a - b)/2, (3b - a)/2\right]\). To recover a single-valued solution, we insert a shock discontinuity.  Due to the form of transformation~\Cref{equation:Transformation}, the flux potential \(\Phi(u)\) always conserved across the shock so that in principle, a shock could be inserted anywhere in the multi-valued region. \arev As a result of the symmetry \(\Phi(u)_x\) is always conserved across the shock. \endrev To specify a second condition, we notice that equation \eqref{equation:Reaction_diffusion} can be rewritten in the form
\begin{equation}
    \label{eq:Phi}%
    \pd{u}{t} = \pdn{}{x}{2}\Phi(u)+R(u),
 \end{equation}
and consequently we impose that \(\Phi'(u) = D(u)\) is also continuous across the shock. That is, the shock position is determined by the following two conditions
\begin{equation}
    \Phi(\ShockMin) = \Phi(\ShockMax), \quad \text{ and } \quad \Phi'(\ShockMin)=D(\ShockMin) = D(\ShockMax)=\Phi'(\ShockMax), 
    \label{equation:Equal_diffusion_condition}
\end{equation}            
where \(\ShockMin\) and \(\ShockMax\) are the left (lower) and right (upper) endpoints of the shock. For quadratic diffusivity~\cref{equation:Diffusion}, these points are
\begin{equation}
    \label{ShockMinMax}%
    \ShockMin = \frac{a + b - \sqrt{3(a - b)^2}}{2} \quad \text{ and } \quad \ShockMax = \frac{a + b + \sqrt{3(a - b)^2}}{2},
\end{equation}
where we also assume that \(b < a(2 + \sqrt{3})\) so that $\ShockMin>0$. Due to the relationship between the diffusivity and the reaction term~\Cref{equation:Reaction}, \(R(u)\) is also preserved across the shock. In addition, both \(u_x\) and \(u_t\) are preserved across the shock. \Cref{figures:RDP_shock} shows the behaviour \arev of \endrev both the diffusivity and the reaction across the shock. In~\Cref{figures:RDP_shock_a}, the dashed horizontal line shows that the solution jumps over both the zeros of the diffusivity, \(u = a, b\) and the region where \(D(u) < 0\) demonstrating that \(\Phi'(u) = D(u)\) is conserved. Likewise, in~\Cref{figures:RDP_shock_b} the solution jumps over the locations where the reaction term is singular, and the region between the vertical lines where the reaction term is positive. \Cref{figures:RDP_shock_c} shows the flux potential around the shock, demonstrating that \(\Phi(u)\) is conserved across the shock.
%
%%%%%%%%%%%%%%%%%%%%%%%%%%%%%%%%%%%%
%%%%% Time-dependent solutions %%%%%
%%%%%%%%%%%%%%%%%%%%%%%%%%%%%%%%%%%%
%
\subsection{Examples of analytic time-dependent solutions}\label{subsection:Example_time_dependent_solutions}      
\begin{figure}[htbp!]
    \centering
    \subcaptionbox{\label{fig:time-dependent_receding}}{\includegraphics[width=0.45\linewidth]{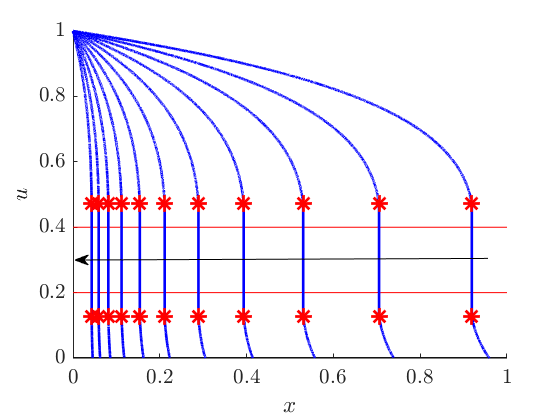}}
    \subcaptionbox{\label{fig:time-dependent_colliding}}{\includegraphics[width=0.45\linewidth]{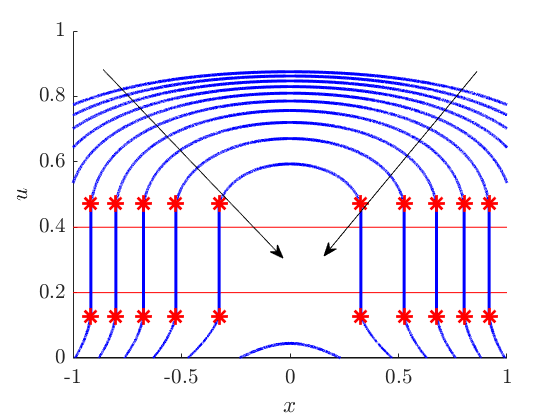}}
    \caption{Time-dependent solutions. Curves (blue) show the solution at different times and the arrows indicate how the solutions move with time. Stars (red) indicate the beginning and end of the vertical shock. Horizontal lines (red) show when the diffusivity is zero. In both solutions \(\kappa = -1\), \(A = 0.08\), \(u^* = 1\), \(a = 0.2\), \(b = 0.4\). (a) Solution which satisfies the condition that at \(u(1,t) = 1\), it has \(c_1 = \Phi(0)\) and \(c_2 = - c_1 = -\Phi(0)\). (b) Colliding wave solution, with \(c_1 = c_2 = \Phi(0)\).} 
    \label{figures:Solution_time_dependent}
\end{figure}

We obtain different analytic solutions by adjusting the constants \(c_1\) and \(c_2\) in~\cref{equation:Implicit_solution}. For example, setting \(c_1 < 0\) and \(c_2 = -c_1\) gives rise to a receding solution that satisfies \(u(0,t) = 1\) for all \(t\). Although the receding solution \(u(x,t)\) gets steeper as the front approaches \(x = 0\), the shock persists because the underlying analytic solution remains multi-valued. This solution is shown in~\Cref{fig:time-dependent_receding}.    Alternatively, if we remove the relationship between \(c_1\) and \(c_2\) and instead merely impose \(c_1 < 0\) and \(c_2 < 0,\) we obtain the colliding wave solution shown in~\Cref{fig:time-dependent_colliding}. In the colliding wave case, there is a time where the solution lies completely below the values of $u$ for which $D(u)<0$ and becomes smooth; at later times the solution becomes zero everywhere.

\brev Given that the density \(u(x,t)\) typically represents a biological or chemical concentration, we usually require \(u(x,t)\ge0\) in feasible solutions to reaction--diffusion equations. The nonclassical symmetry solutions presented in~\Cref{figures:Solution_time_dependent} do not have \(u \geq 0\) for all values of \(x.\) However, we can obtain solutions with non-negative \(u\) by restricting the domain \cite{Compactly_Supported_Solutions:edwards2018}. For example, the time-dependent receding solution in~\Cref{fig:time-dependent_receding} has non-negative density for \(x \in [0, L(t)],\) where \(x = L(t)\) is a moving boundary. When viewed on these restricted domains, the solutions are sharp-fronted. That is, at the boundary \(x = L(t)\) where the density \(u(L(t), t) = 0,\) the gradient is non-zero, \(u_x(L(t), t) \neq 0.\) The region on which \(u\) is non-negative evolves with time, such that the function \(L(t)\) describes the position of the sharp front. An important comment is that we did not impose a condition on \(u_x\) at the moving boundary at the outset of the nonclassical symmetry analysis. Instead, we infer the condition from our analytic solution (see below), and it arises as a consequence of the non-negative \(u\) requirement. \endrev

\brev Viewing the solution to the PDE as one obtained from a moving-boundary problem would involve replacing a far-field boundary condition with a condition at \(x = L(t).\) \endrev A well-known moving-boundary problem \brev in this spirit \endrev is the classical Stefan problem. \brev This problem \endrev involves solving the linear heat equation on \(x \in [0, L(t)],\) where the boundary evolves according to the Stefan condition~\Cref{equation:StefanCondition}                
\begin{equation}        
    \left.\fd{u}{x}\right\rvert_{x = L(t)} = -\beta\fd{L}{t},
    \label{equation:StefanCondition}
\end{equation}    
where \(\beta\) is the Stefan parameter. We can calculate both the speed of the moving boundary and the flux at the moving boundary exactly \brev using our analytic solutions. Doing so reveals that \endrev our solutions do not satisfy the regular Stefan condition~\eqref{equation:StefanCondition}, but instead satisfy an alternate condition whereby the speed is inversely proportional to the gradient,
\begin{equation}        
    \left.\pd{u}{x}\right\rvert_{x = L(t)} = \frac{\kappa\Phi(0)}{L'(t)}. \label{equation:StefanCondition_Specific}
\end{equation}
Both the time-dependent and colliding wave solutions in \Cref{figures:Solution_time_dependent} satisfy this condition (in the case of colliding waves, both boundaries satisfy the condition). \brev Although the physical or biological interpretation of this Stefan-like condition~\cref{equation:StefanCondition_Specific} remains to be determined, moving boundary conditions with boundary speed proportional to the reciprocal of the density gradient have been reported in the literature~\citep{Lundberg2013}. \endrev 

\brev One consequence of our Stefan-like condition~\cref{equation:StefanCondition_Specific} is that \(L'(t) < 0\) in travelling wave and time-dependent solution. This is because the gradient \(u'(x)\) is negative at \(x = L(t)\) (see~\Cref{fig:time-dependent_receding,figures:Travelling_wave_solution_nozoom}), and \(\kappa\) and \(\Phi(0)\) are both negative constants so that \(\kappa \Phi(0)\) is positive. This gives rise to receding density profiles. Although receding fronts are uncommon in reaction--diffusion models, the receding behaviour is similar to that seen in recent works that adapt the classical Stefan problem to solve reaction--diffusion (Fisher-KPP) equations on moving boundaries~\cite{Spreading-vanishing:du2010, Revisiting_the_Fisher:el-hachem2019, Invading_and_receding:el-hachem2021,Semi-Infinite:fadai2021}. Due to \endrev the inverse relationship between the gradient and the speed of the boundary~\cref{equation:StefanCondition_Specific}, steeper density gradients lead to slow-moving fronts, whereas shallower gradients increase front speed. Large flux at the boundary \((D(u)u_x)|_{x = L(t)}\) thus corresponds to slow fronts, and small flux corresponds to fast fronts.
%
%%%%%%%%%%%%%%%%%%%%%%%%%%%%%%%%%%%%%%%%%%%%%%
%%%%% Receding travelling wave solutions %%%%%
%%%%%%%%%%%%%%%%%%%%%%%%%%%%%%%%%%%%%%%%%%%%%%
%
\subsection{Analysis of a receding travelling wave solution} \label{subsection:An_example_TW}
\begin{figure}[htbp!]
    \centering
    \subcaptionbox{\label{figures:Travelling_wave_solution_nozoom}}{\includegraphics[width=0.45\linewidth]{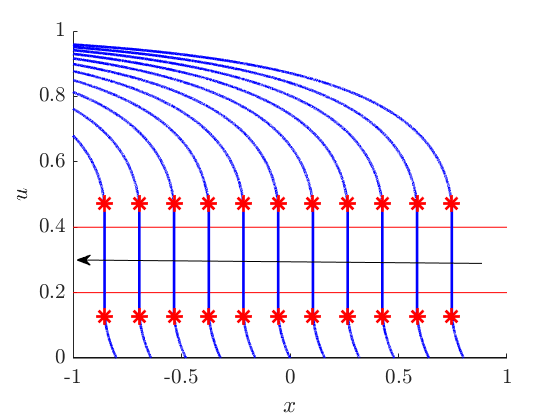}}
    \subcaptionbox{\label{figures:Travelling_wave_solution_zoom}}{\includegraphics[width=0.45\linewidth]{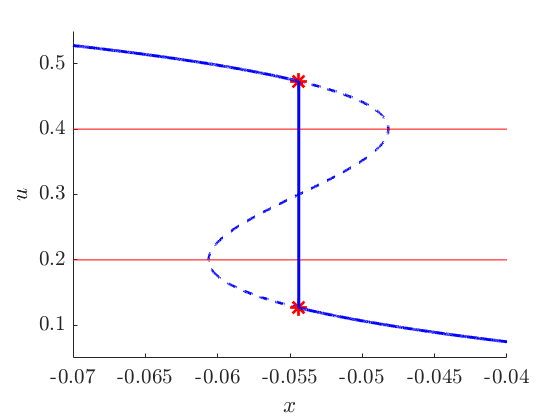}}
    \caption{Travelling wave solution. Curves (blue) show the solution at different times and the arrow indicates how the solution moves with time. Stars (red) indicated the beginning and end of the vertical shock. Horizontal lines (red) show where the diffusivity is zero. (b) The dashed curve shows the multi-valued part of the solution at a single time. Here \(\kappa = -1\), \(A = 0.08\), \(u^* = 1\), \(a = 0.2\), \(b = 0.4\), \(c_1 = \Phi(0)\) and \(c_2 = 0\).}
    \label{figures:Travelling_wave_solution}
\end{figure}

While the time-dependent solutions in~\Cref{section:ExampleSolution} involve non-zero \(c_1\) and \(c_2,\) we can obtain a travelling wave solution by choosing \(c_1 < 0\), \(c_2 = 0.\) The implicit solution~\Cref{equation:Implicit_solution} for \(\Phi(u)\) then becomes
\begin{equation}
    \label{eq:Implicit_Solution_TW}%
    \Phi(u) = c_1\e^{kx + At}.
\end{equation}
We can write Equation~\cref{eq:Implicit_Solution_TW} in terms of the travelling wave coordinate \(z=x-ct\), where \(c=-A/k=-kD(0)\). The travelling wave solution then has an implicit form, expressed in terms of \(z\) as
\begin{equation}
    z = x -ct= \frac{1}{k}\log{\left(\frac{\Phi(u)}{c_1}\right)}.
    \label{equation:implicit_travelling_wave_solution}
\end{equation} 
As described above, the travelling wave solution~\cref{equation:implicit_travelling_wave_solution} obtained using the nonclassical symmetry is multi-valued around the region where the diffusivity takes negative values. \Cref{figures:Travelling_wave_solution} shows a shock inserted so that the diffusion is constant across the shock, that is, according to conditions \eqref{equation:Equal_diffusion_condition}. The start and end of the shock are marked by stars (red) while the horizontal lines (red) show the zeros of the diffusivity; between these lines, the diffusivity is negative. The arrow shows that this travelling wave solution moves leftward (\(c < 0\)), and is a receding front. This is a consequence of requiring \(R(0)=0,\) which means that \(A = -\kappa D(0) > 0.\) The shock avoids both the region of negative diffusivity and the singularities in the reaction term, and both the diffusivity and reaction terms are continuous across the shock. The behaviour of \(D(u)\) and \(R(u)\) across the shock corresponds to the horizontal dashed lines in~\Cref{figures:RDP_shock_a} and~\Cref{figures:RDP_shock_b} respectively.
            
\begin{figure}[htbp!]
    \centering
    \subcaptionbox{Boundary position, \(L(t).\) \label{figures:boundary}}{\includegraphics[width=0.45\linewidth]{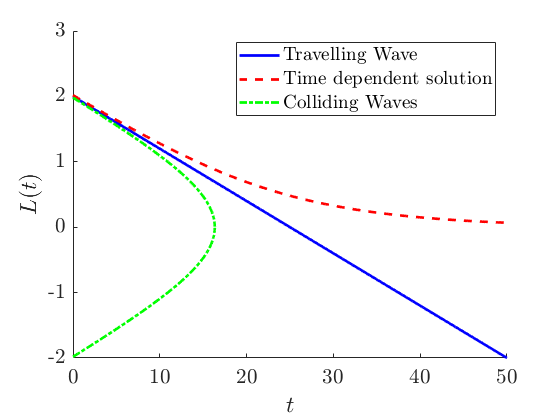}}
    \subcaptionbox{Flux, \(-D(u)u_x\) at \(x = L(t).\) \label{figures:flux}}{\includegraphics[width=0.45\linewidth]{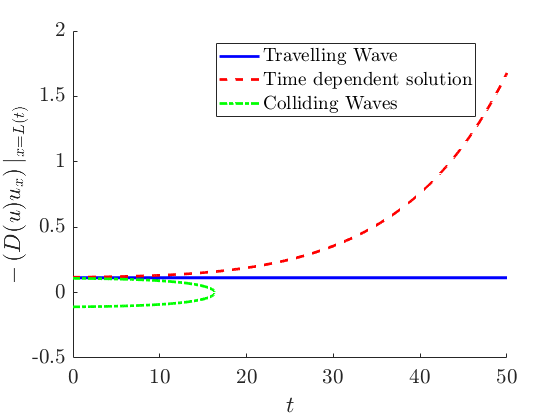}}
    \caption{Comparison of the boundary position \(L(t)\) and the flux \(-D(u)u_x\) at \(x = L(t)\) for the receding time-dependent, colliding wave, and travelling wave solutions. The dash-dot curves (green) are multi-valued because it corresponds to the colliding wave solution and has two boundaries. Solutions obtained using parameter values \(\kappa= - 1\), \(A = 0.08\), \(u^* = 1\), \(a = 0.2\), \(b = 0.4\). For the time-dependent solution, we set \(c_1 = - c_2 = -0.0153\), for the colliding wave solution we have \(c_1 = c_2 = -0.0153\), and for the travelling wave solution we use \(c_1 = -0.0153\), \(c_2 = 0\).}
    \label{figures:boundary_and_flux}
\end{figure}

Like the time-dependent solutions in~\Cref{figures:Solution_time_dependent}, the travelling wave solution in~\Cref{figures:Travelling_wave_solution} also satisfies the Stefan-like condition~\Cref{equation:StefanCondition_Specific} at the location where \(u = 0.\) Choosing \(c_1 = \Phi(0)\) means that the sharp front in \(u\) coincides with \(z = 0\) (that is, \(u(0) = 0\)). \Cref{figures:boundary} shows the location of the boundary and \Cref{figures:flux} shows the flux at the boundary. For the travelling wave solution, the speed at the moving boundary and the flux at the boundary are both constant. In contrast, time-dependent solutions have an inverse relationship between the speed of the moving boundary and the flux. In the receding time-dependent solution, as the moving boundary approaches zero its speed decreases while the flux grows. For the colliding wave solution, the dash-dot curves (green) in \Cref{figures:boundary_and_flux} are multi-valued because of the two contact lines; as the two boundaries approach zero the speed of the moving boundaries starts to increase and the flux at both boundaries approaches zero. 

The time-dependent and travelling wave solutions are both left-moving (receding) solutions. \citet{Shock-fronted:li2021} derive a necessary condition for the existence of left-moving travelling waves satisfying the usual far-field boundary conditions (that is, \(u\to 0,1\) and \(u_z\to 0\) as \(z\to\pm\infty\) respectively). For a sharp-fronted wave, the condition of~\citet{Shock-fronted:li2021} changes from
\begin{equation}
    \int_0^{u_a} D(u)R(u) \df{u} < -\frac{1}{2}\left(g(u_a)\right)^2 < 0, \quad \text{for all } u_a \in (0, \ShockMin), 
    \label{equation:left_condition}
\end{equation}
to
\begin{equation}
    \int_0^{u_a} D(u)R(u) \df{u} < -\frac{1}{2}\left(g(u_a)\right)^2 + \frac{1}{2}\left(g(0)\right)^2, \quad \text{for all } u_a\in(0,u_l), 
    \label{equation:left_condition_our_version}
\end{equation}
where \(g(u)=D(u)u_z\). Consequently, while a model with a reaction term that is positive in $(0,u_a)$ cannot support a smooth-fronted receding travelling wave (that is, a solution where  \(u_z \to 0\) as \(z\to\infty\)), a compactly-supported sharp-fronted receding travelling wave (for which $g(0)\ne0$) can be supported.
            
To analyse the receding travelling wave solution further, we apply the transformation \(z = x - ct\) to the RDE~\cref{equation:Reaction_diffusion} to obtain the ordinary differential equation
\begin{equation}
    \label{equation:Travelling_wave_ode}%
    \fd{}{z}\left(D(u)\fd{u}{z}\right) + c \fd{u}{z} + R(u) = 0.
\end{equation}
The implicit solution~\cref{equation:implicit_travelling_wave_solution} with \(c_2 = 0\) is a travelling wave solution \(u(z)\) that satisfies~\cref{equation:Travelling_wave_ode}. Following~\citet{Travelling_wave_solutions:li2020}, we write~\cref{equation:Travelling_wave_ode} as the system
\begin{subequations}
    \label{equation:phase_plane}%
    \begin{align}
        D(u) \fd{u}{z} &= q, \label{equation:phase_plane_u}\\
        D(u) \fd{q}{z} &= -cq - D(u)R(u). \label{equation:phase_plane_p}
    \end{align}
\end{subequations}
The left-hand sides of~\cref{equation:phase_plane} vanish when \(D(u) = 0,\) giving rise to a {\it wall of singularities} at \(u = a\) and \(u = b\) \citep{Existence_of_traveling_wave:harley2014,Travelling_wave_solutions:li2020,Lokta-Volterra_equations:pettet2000,Folds_canards:wechselberger2010}, \emph{i.e.}, solutions of ~\cref{equation:phase_plane} will, in general, cease to exist when reaching these walls of singularities due to a finite-time blow-up. Possible exceptions are points known as {\it holes in the walls}~\citep{Existence_of_traveling_wave:harley2014,Travelling_wave_solutions:li2020,Lokta-Volterra_equations:pettet2000,Folds_canards:wechselberger2010} where the left-hand and right-hand sides of~\cref{equation:phase_plane} both vanish. However, the conditions $D(u)=0$ and $q=0$ do not correspond to holes in the wall; although the right-hand side of~\cref{equation:phase_plane_u} vanishes when $q=0$, the right-hand side of~\cref{equation:phase_plane_p} does not vanish when $D(u)=0$ since the product
\[D(u)R(u) = \left(A + \kappa D(u)\right)\int_{u^*}^u D(u') \df{u'}=\left(A + \kappa D(u)\right)\Phi(u),\]
is non-zero at the walls of singularities, \emph{i.e.}, $A\int_{u^*}^{u=a,b} D(u') \df{u'} \neq 0$. Thus, smooth solutions between $u=0$ and $u=1$ are not possible in system~\cref{equation:phase_plane} and shocks are a necessary ingredient.
            
\begin{figure}[htbp!]
    \centering
    \includegraphics[width=0.7\linewidth]{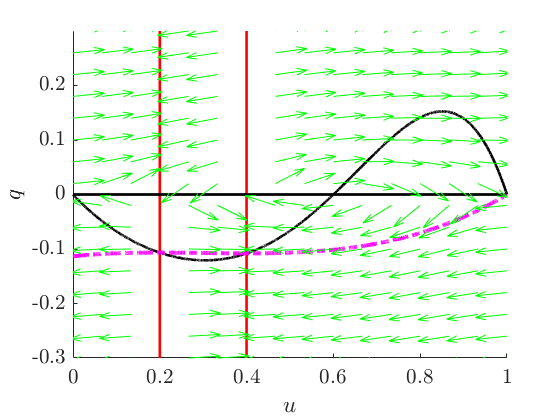}         
    \caption{Phase portrait of \Cref{equation:phase_plane}. The arrows (green) indicate the direction field of the system. The dashed curve (magenta) is the solution trajectory of our analytic solution, which emanates from \((1,0)\) and terminates at \((0,-q_0)\) satisfying the boundary condition~\cref{equation:StefanCondition_Specific}. The vertical lines (red) are walls of singularities. The solid horizontal line and the solid curve (black) are the nullclines of the system.}
    \label{figure:phase_plane}
\end{figure}

Nevertheless, we can use our analytic multi-valued solution to calculate solution trajectories in the phase space. \Cref{figure:phase_plane} shows a phase plane portrait of \Cref{equation:phase_plane} with the trajectory of our multi-valued solution. The phase plane also shows why a shock is required: between the walls the diffusivity is negative and the direction of the flow switches. As a consequence, our solution cannot be represented by a single trajectory and instead is a combination of two solution branches. The phase portrait in~\Cref{figure:phase_plane} also demonstrates that the sharp-fronted receding travelling wave has non-zero flux at the moving boundary (where \(u = 0\)).
%
%%%%%%%%%%%%%%%%%%%%%%%%%%%%%%
%%%%% Stability analysis %%%%%
%%%%%%%%%%%%%%%%%%%%%%%%%%%%%%
%
\subsection{Stability analysis}
We analyse the spectral stability of the constant and travelling wave solutions to \cref{equation:Reaction_diffusion}, relative to perturbations in an appropriately chosen space. To that end, we first derive the {\em linearised equation}, linearised about a solution $\bar{u}$, which we write in the form $p_t = \cL(\bar{u}) p$, where $p$ is a perturbation to the solution $\bar{u}$ and $\cL(\bar{u})$ is the {\em linearised operator}. To derive the linearised equation, we consider solutions to \cref{equation:Reaction_diffusion} of the form $u = \bar{u} + \ve p(x,t),$ where \(\ve \ll 1.\) Collecting terms of $\mathcal{O}(\ve)$ as \(\ve \to 0,\) we arrive at a PDE for the perturbation $p$, 
\begin{equation}
    \label{equation:linearised}%
    p_t = \left(D(\bar{u})p\right)_{xx} + R'(\bar{u})p = : \cL(\bar{u}) p.
\end{equation}
For a general solution of \cref{equation:Reaction_diffusion}, the right hand side of \cref{equation:linearised}, and in particular the operator $\cL(\bar{u})$ will depend on both $x$ and $t$. However, for the special solutions we consider, $\cL$ will either be a constant coefficient operator, or depend on a single variable only. 
\subsubsection{Linear stability analysis of the constant solutions}\label{subsubsection:Stabillity_constant} 
We suppose first that our solution $\bar{u}$ is a constant solution, either $0$ or $1$. Then, \cref{equation:linearised} becomes a constant-coefficient second-order PDE, which is solvable via separation of variables. Supposing that $p(x,t)$ takes the form $p(x,t)= e^{\lambda t} p(x)$ we obtain the ordinary differential equation
\begin{equation}
    \label{equation:linear_const}
    \lambda p =  D(\bar{u})p_{xx} + R'(\bar{u})p.
\end{equation}
We choose the domain of our perturbations to be $p \in L^2(\R)$, which provides the boundary conditions $\lim_{x \to \pm \infty} p = 0$, and enables us to solve \cref{equation:linear_const} via the Fourier transform. We are led to the dispersion relation, 
\begin{equation}
    \lambda = -\alpha^2 D(\bar{u}) + R'(\bar{u}), 
    \label{equation:constant_lambda_condition_0}
\end{equation}
where $\alpha$ is the Fourier variable in space. If $\Re(\lambda) < 0$ then, to first order, the perturbation will decay to the constant solution $\bar{u}$, indicating that \(u=0\), or $u=1$ is linearly stable to $L^2(\R)$ perturbations. In contrast, if $\Re(\lambda) > 0$ the perturbation to the constant solution \(\bar{u} \) will grow with time. 

For $u = 0$, the quadratic diffusivity satisfies \(D(0) > 0,\) so $-\alpha^2 D(0) < 0$ for all real wave numbers \(\alpha.\) Stability of the constant solution $u = 0$ then depends on the sign of $R'(0)$. The derivative of the reaction term (with $A = -\kappa D(0)$) is
\begin{equation}
    R'(u) = \left(\frac{\kappa D(0)D'(u)}{D(u)^2}\right) \Phi(u) + \left(\frac{-\kappa D(0)}{D(u)} + \kappa\right) D(u). 
    \label{equation:reaction_derivative}
\end{equation}
At $u = 0$ we have
\begin{equation}
    R'(0) = \left(\frac{\kappa D'(0)}{D(0)}\right) \Phi(0). \label{equation:reaction_derivative_at_zero}
\end{equation}
Since $\kappa < 0$, $\Phi(0) < 0$, and $D(0) > 0$, the derivative $R'(0)$ has the same sign as $D'(0) = - a - b$, which is always negative as $0 < a < b$. Therefore, small perturbations around the constant solution $u = 0$ will decay, and the constant solution $u = 0$ is (linearly) stable.
                
Likewise, we can consider the stability of the constant solution $u = 1$. Since $D(1) > 0$, stability depends on the value of $R'(1)$, which is
\begin{equation}
    R'(1) = -\kappa D(1)\left(\frac{D(0)}{D(1)} - 1\right). 
    \label{equation:reaction_derivative_at_one}
\end{equation}
Since $-\kappa D(1) > 0$ for quadratic diffusivity, the sign of $R'(1)$ depends on the sign of $\left(D(0)/D(1) - 1\right).$ This term is positive if $D(1) < D(0)$, negative if $D(0) < D(1)$, and is zero if $D(1) = D(0)$. Therefore, $R'(1)$ is positive if $a + b > 1$, is negative if $a + b < 1$ and is zero when $a + b = 1$. Quadratic diffusivity with \(a + b > 1\) thus admits long-wavelength instabilities for sufficiently small \(\alpha.\) The solutions presented in \Cref{figures:Travelling_wave_solution} are for $a+b<1$, for which so $u=1$ is stable relative to perturbations in $L^2(\R)$.
\subsubsection{Stability analysis of the receding travelling wave solution}\label{subsubsection:TWS_Stability}
Having considered constant \(\Bar{u}\), we now address the linearisation of \cref{equation:Reaction_diffusion} in the moving frame $(z,t) = (x-c t, t)$, about the receding travelling wave solution described in \Cref{subsection:An_example_TW}. Again, we subject the solution to perturbations of the form $p(z)e^{\lambda t}$, where $p(z)$ will be in an appropriate space. Our equation for the perturbation is
\begin{equation}
    \lambda p = \cL(\bar{u})p := \left(D(\bar{u})p\right)_{zz} + cp_z + R'(\bar{u})p.
    \label{equation:travelling_wave_perturbation_ode}
\end{equation}
The right hand side of \cref{equation:travelling_wave_perturbation_ode} depends only on \(z\). We require that $p(z)$ be a solution to \cref{equation:travelling_wave_perturbation_ode} which is square integrable on $(-\infty,0)$ which also vanishes at $z = 0$ (\emph{i.e.} in $L_0^2((-\infty,0))$). This means that we impose the boundary conditions $p(0) = 0$ and $\lim_{z\to -\infty} p(z) =0$. 

To uniquely solve the linearised ODE \cref{equation:travelling_wave_perturbation_ode}, we need to impose conditions across the shock. To this end, we ask that $p$ be $C^1$ across the shock. We note that this means that the coefficient of the highest derivative of $D(\bar{u})$ will be continuous across the shock, but the lower order coefficients may not be. We remark that since the shock occurs at a single point, the set of functions which are $C^1$ across the shock is still a densely defined subset of $L_0^2((-\infty,0))$, so imposing this condition will not affect the spectrum of $\cL(\bar{u})$. 

For our discussion on spectral stability of the receding travelling wave solution, we use the language and notation outlined in Chapters 2 and 3 of the textbook by \citet{Spectral_and_dynamical:kapitula2013}. The set of $\lambda \in \C$ such that $\cL(\bar{u}) - \lambda$ does not have a bounded inverse on $L_0^2((-\infty,0))$, will be called {\em the spectrum} of $\cL(\bar{u})$, and denoted $\sigma(\cL(\bar{u}))$. We note that the spectrum of a closed linear operator on a Hilbert space is a (topologically) closed set, so if we can find an open set which is in the spectrum, then its closure must be as well. We will call the wave $\bar{u}$ {\em spectrally stable}, if $\lambda \in \sigma(\cL(\bar{u}))$ means that $\Re(\lambda) < 0$. 
\begin{theorem} 
    \label{th:stability} 
    The travelling wave solution found in~\Cref{subsection:An_example_TW} is spectrally stable, provided \(a\) and \(b\) are such that  \(\frac{1}{2} D(\bar{u})_{zz} + R'(\bar{u}) < 0\) for \(u \in [0,\ShockMin]\cup[\ShockMax,1].\)
\end{theorem}
\begin{remark} 
    Due to the condition \(\frac{1}{2} D(\bar{u})_{zz} + R'(\bar{u}) < 0,\) we require \(R'(1) < 0\).
\end{remark}
\begin{remark} 
    Sometimes waves are referred to as spectrally stable if they have spectrum in the left-half plane, {\em and} perhaps a simple eigenvalue associated with translation invariance at the origin. As our perturbations are in $L^2_0((-\infty,0))$, the eigenvalue associated with translation invariance $\lambda = 0$ will not be an eigenvalue of the linearised operator. 
\end{remark}
\begin{proof}
    We begin with an analysis of the so-called far-field operator of \cref{equation:travelling_wave_perturbation_ode}. Letting $z \to -\infty$, we have that solutions to \cref{equation:travelling_wave_perturbation_ode} will asymptotically behave like solutions to the constant coefficient ordinary differential equation
    \begin{equation}
        \lambda p = D(1) p_{zz} - k D(0) p_z + R'(1) p,
        \label{equation:travelling_wave_perturbation_u=1}
    \end{equation}
    where we recall that $c = - k D(0)$. The characteristic equation of ODE \cref{equation:travelling_wave_perturbation_u=1} is 
    \begin{equation}
        \lambda = D(1) r^2 - k D(0) r + R'(1),
    \end{equation}
    with characteristic exponents given as the roots
    \[r_\pm := \frac{ k D(0) \pm \sqrt{ k^2 D(0)^2 - 4 D(1) (R'(1) - \lambda)} }{2 D(1)}.\]
    We are interested in the sign of the real parts of $r_\pm$ for a given value of $\lambda$. Supposing that $r = i \alpha$ for $\alpha \in \R$, we have the dispersion relation
    \begin{equation} 
        \label{equation:dispersion_relation}%
        \lambda = -\alpha^2 D(1) - i kD(0) \alpha + R'(1).
    \end{equation}
    Equation~\Cref{equation:dispersion_relation} is a parametrised curve for $\lambda$ in the complex plane, parametrised by $\alpha$, and divides the complex plane into two disjoint regions. We denote the region to the right of the curve by $\Omega$, while the region to the left of the curve we will call $B$, see~\Cref{figure:spectral_stability}. For $\lambda \in \Omega$ we have that $\Re(r_+) >0$, while $\Re(r_-) <0$. However, for $\lambda \in B$ we have that both $\Re(r_\pm) >0$, meaning that $\cL(\bar{u})-\lambda$ is not invertible here. As both $r_\pm$ have positive real parts, we have a spanning set of solutions to \cref{equation:travelling_wave_perturbation_ode} which decay to $0$ as $z \to -\infty$. Then the solution satisfying the right-hand boundary condition would necessarily be a linear combination of these solutions, and would hence lie in the kernel of $\cL(\bar{u}) - \lambda$. Owing to our remark above about the spectrum being closed, we conclude that the spectrum at least contains the closure of the set $B$. 

    We have that \(D(1) > 0 \), so the maximum real value of \(\lambda\) in \cref{equation:dispersion_relation} depends on \(R'(1)\), which we know is negative provided that the conditions in the statement of~\Cref{th:stability} are satisfied. Thus this part of the spectrum at least is contained entirely in the left half plane. 
    \begin{figure}[htbp!]
        \centering
        \includegraphics[width=0.75\linewidth]{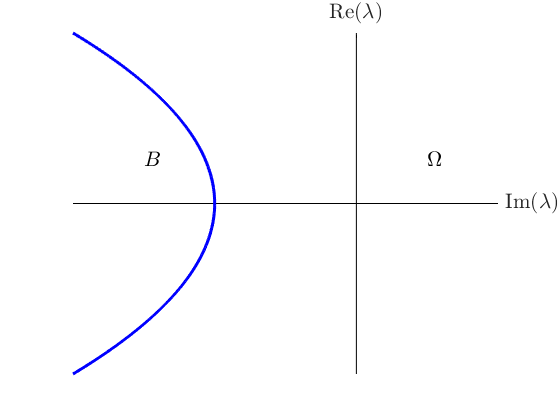}
        \caption{Schematic of the dispersion relation~\cref{equation:dispersion_relation} in the complex plane dividing the plane into two disjoint regions, denoted \(B\) and \(\Omega.\)}
        \label{figure:spectral_stability}
    \end{figure}
    \begin{remark} 
        In fact, what the above calculation partially shows is that the so-called {\em Fredholm index} of $\cL(\bar{u}) - \lambda$ is not zero for $\lambda \in B$ (it is in fact 2), but as we can show lack of invertibility of $\cL(\bar{u})- \lambda$ for $\lambda \in B$ directly, we avoid the diversion into functional analysis terminology. For a more complete description of why the previous calculation computes the Fredholm index of $\cL(\bar{u})- \lambda$, we refer the reader to \cite{Spectral_and_dynamical:kapitula2013}. 
    \end{remark}
    
    The previous computation shows that for $\lambda \in \Omega$ we can compute the Green's function of \cref{equation:travelling_wave_perturbation_ode}, and thus invert $\cL(\bar{u}) - \lambda$, provided we do not have an eigenvalue/eigenfunction to \cref{equation:travelling_wave_perturbation_ode}. We first show that any eigenvalues must be real, and then show that for any \(a\) and \(b\) satisfying the conditions in the statement of~\Cref{th:stability} the eigenvalues must be negative.                                       

    To see that eigenvalues of \cref{equation:travelling_wave_perturbation_ode} with $\lambda \in \Omega$ must be real, we multiply \cref{equation:travelling_wave_perturbation_ode} by the integrating factor 
    \begin{equation}
        D(\bar{u})^2\exp{\left(- k D(0)\int \frac{1}{D(\bar{u})}\, dz\right)},
    \end{equation}
    so that equation~\Cref{equation:travelling_wave_perturbation_ode} is put into Sturm--Liouville form
    \begin{equation}
        \label{eq:Sturm}
        \cM p := (S_1(z)p')' + S_2(z)p = \lambda S_3(z) p,
    \end{equation}
    where the three \(z\)-dependent functions are
    \begin{subequations}
        \begin{align}
            S_1 &= D(\bar{u})^2\exp{\left(- k D(0)\int \frac{1}{D(\bar{u})}\, dz\right)} \\
            S_2 &= D(\bar{u})\left(R'(\bar{u}) + D(\bar{u})_{zz}\right)\exp{\left(- k D(0)\int \frac{1}{D(\bar{u})}\, dz\right)} \\
            S_3 &= D(\bar{u})\exp{\left(- k D(0)\int \frac{1}{D(\bar{u})}\, dz\right)}.
        \end{align}
    \end{subequations}
    Noting that $S_3$ is a positive, bounded, continuous function (even across the shock because $u$ always takes values where $D(u)>0$ and $D(u)$ is continuous across the shock~\cref{equation:Equal_diffusion_condition}), a calculation shows that for $\lambda \in \Omega$, if $\lambda$ is an eigenvalue of $\cL(\bar{u})$, then it will be an eigenvalue of $\frac{1}{S_3(z)} \cM$ in the weighted Hilbert space $(L^2_0((-\infty,0)))_{S_3}$ with the (usual) weighted norm
    \[ <v, v>_{S_3} := \int v \bar{v} S_3 \df{z}. \]
    Indeed we have that for $z \ll -1$, $D(\bar{u}) \sim D(1)$. Solutions when $\lambda \in \Omega$ will decay like $e^{\Re{(r_+)z}}$, and we have that $2 \Re{(r_+)} - k D(0)/D(1) >0$ when $\lambda \in \Omega$. Therefore, eigenfunctions when $\lambda \in \Omega$ will be in $(L^2_0((-\infty,0)))_{S_3}$ (note: this does not happen for $\lambda \in B$ or its closure). But such an operator is self-adjoint in this space, and hence only has real eigenvalues. Thus for $\lambda \in \Omega$, if it is an eigenvalue of $\cL(\bar{u})$ in $L^2_0((-\infty,0))$, then it must be real. 
                
    Lastly to see that any eigenvalues $\lambda \in \Omega$ must be negative, we note that if $\lambda$ is real, the corresponding eigenfunction will also be real.  We multiply \Cref{equation:travelling_wave_perturbation_ode} by $p$ and integrate,
    \begin{align}
                        \int_{-\infty}^{z_s} \lambda p^2 \df{z} + \int_{z_s}^0 \lambda p^2 \df{z} 
                        =& \int_{-\infty}^{z_s} \left(D(\bar{u})p\right)_{zz} p \df{z} + \int_{z_s}^0 \left(D(\bar{u})p\right)_{zz} p \df{z} \nonumber \\
                        &+ \int_{-\infty}^{z_s} c p p_z \df{z} + \int_{z_s}^0 c p p_z \df{z} 
                        + \int_{-\infty}^{z_s} R'(\bar{u})p^2 \df{z} + \int_{z_s}^0 R'(\bar{u})p^2 \df{z}, \label{equation:test_integral_intial}
    \end{align}
    where $z_s$ is the location of the shock. Integrating the speed terms results in
    \begin{equation}
        \int_{-\infty}^{z_s} c p p_z \df{z} + \int_{z_s}^{0} c p p_z \df{z} = \frac{c}{2}\left(p(z_s^-)^2 - p(z_s^+)^2\right), 
        \label{equation:test_integral_speed}
    \end{equation}
    where \(p(z_s^-)\) is the value of \(p(z)\) as \(z \rightarrow z_s\) from below and \(p(z_s^+)\) is the value of \(p(z)\) as \(z \rightarrow z_s\) from above. For the diffusion term, using integration by parts twice we arrive at
    \begin{align}
        \int_{-\infty}^{z_s} \left(D(\bar{u})p\right)_{zz} p \df{z} + \int_{z_s}^0 \left(D(\bar{u})p\right)_{zz} p \df{z} =& \int_{-\infty}^{z_s} \frac{1}{2} D(\bar{u})_{zz} p^2 \df{z} + \int_{z_s}^0 \frac{1}{2} D(\bar{u})_{zz} p^2 \df{z} \nonumber \\
        & -\int_{-\infty}^{z_s} D(\bar{u})p_z^2 \df{z} -\int_{z_s}^0 D(\bar{u})p_z^2 \df{z} \nonumber \\
        & +\left[\left(D(\bar{u})p\right)_{z} p\right]_{z = z_s^{-}} - \left[\left(D(\bar{u})p\right)_{z} p\right]_{z = z_s^{+}} \nonumber \\                    
        & -\left[\frac{1}{2}D(\bar{u})_z p^2\right]_{z = z_s^{-}} + \left[\frac{1}{2}D(\bar{u})_z p^2\right]_{z = z_s^{+}}. \label{equation:test_integral_diffusion}
    \end{align}
    In total we have
    \begin{align}
        \int_{-\infty}^{z_s} \lambda p^2 \df{z} + \int_{z_s}^0 \lambda p^2 \df{z} =& \int_{-\infty}^{z_s} \left(\frac{1}{2} D(\bar{u})_{zz} + R'(\bar{u})\right) p^2 \df{z} + \int_{z_s}^0 \left(\frac{1}{2} D(\bar{u})_{zz} + R'(\bar{u})\right) p^2 \df{z} \nonumber \\
        & - \int_{-\infty}^{z_s} D(\bar{u})p_z^2 \df{z} - \int_{z_s}^0 D(\bar{u})p_z^2 \df{z} + S, \label{equation:test_integral_final}
    \end{align}                    
    where
    \begin{align}
        S =&\, \frac{c}{2}\left(p(z_s^-)^2 - p(z_s^+)^2\right) + \left[\frac{1}{2} D(\bar{u})_z p^2\right]_{z = z_s^{-}} - \left[\frac{1}{2} D(\bar{u})_z p^2\right]_{z = z_s^{+}} \nonumber \\
        & + \left[D(\bar{u}) p p_z\right]_{z = z_s^{-}} - \left[D(\bar{u})p p_z\right]_{z = z_s^{+}}. \label{equation:test_intergral_remainder}
    \end{align}
    We have \(D(\bar{u})_z = D'(\bar{u})\bar{u}_z\), and for quadratic diffusion \(D'(\bar{u}) > 0\) when \(\bar{u} > (a + b)/2\) which is in the region to the left of the shock \(z \in (-\infty, z_s)\) and \(D'(\bar{u}) < 0\) when \(\bar{u} < (a + b)/2\) which is in the region to the right of the shock \(z \in (z_s, 0]\). Also \(\bar{u}_z <0\), therefore \(\left[D(\bar{u})_z p^2/2\right]_{z = z_s^{-}}\) and \(- \left[D(\bar{u})_z p^2/2\right]_{z = z_s^{+}}\) are both negative. The term \(\left(p(z_s^-)^2 - p(z_s^+)^2\right)\) will vanish if \(p\) is continuous at the shock, and because \(D(\bar{u})\) is conserved across the shock the term \(\left[D(\bar{u}) p p_z\right]_{z = z_s^{-}} - \left[D(\bar{u})p p_z\right]_{z = z_s^{+}}\) will vanish if \(p_z\) is also continuous across the shock. Alternatively this entire term will vanish if \(p\) also vanishes at the shock.

    With \(S < 0\) the right-hand side of Equation~\Cref{equation:test_integral_final} will always be negative since \(D(\Bar{u})\) is always positive (since the shock avoids the regions of negative diffusivity) and \(D(\Bar{u})_{zz}/2 + R'(\Bar{u})\) is assumed to be negative. We then need to have \(\lambda < 0\) for the left-hand side of~\cref{equation:test_integral_final} to be negative as well. Hence, any real eigenvalues of $\cL(\bar{u})$ in $\Omega$ will be negative. This completes the proof of~\Cref{th:stability}.
\end{proof}

\begin{figure}
    \centering
    \includegraphics{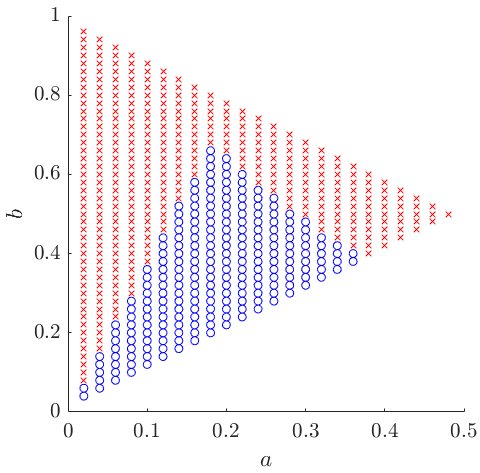}
    \caption{Values of \(a\) and \(b\) for which \(a + b < 1\) and \(0 < a < b < 1\). Circles (blue) indicate the pairs that satisfy the conditions for spectral stability outlined in~\Cref{th:stability}, and the condition \(b < a(2 + \sqrt{3}),\) which guarantees $\ShockMin>0$~\cref{ShockMinMax}. Crosses (red) indicate combinations of \(a\) and \(b\) that do not satisfy both conditions. Our analysis proves spectral stability for \((a,b)\) pairs denoted by blue circles, but makes no statement about the stability of \((a,b)\) pairs denoted by red crosses.}
    \label{figure:allConditions}
\end{figure}

We have proven that receding travelling waves with quadratic diffusivity are spectrally stable whenever \(D(\Bar{u})_{zz}/2 + R'(\Bar{u}) < 0.\) We now explore numerically the values of \(a\) and \(b\) for which the stability analysis applies. \Cref{figure:allConditions} shows \((a,b)\) pairs where \(a + b < 1\) and \(0 < a < b < 1\). The circles (blue) indicate the pairs of \((a,b)\) which satisfy the stability condition \(D(\Bar{u})_{zz}/2 + R'(\Bar{u}) < 0\) as well as \(b < a(2 + \sqrt{3})\) so that the lower shock value \(u_l > 0\) \Cref{ShockMinMax}. The crosses (red) indicate where these conditions are not all satisfied. The left border of the circle (blue) region is determined by the condition \(b < a(2 + \sqrt{3})\),  while the top right border is determined by the stability condition \(D(\Bar{u})_{zz}/2 + R'(\Bar{u}) < 0.\)   

%
%%%%%%%%%%%%%%%%%%%%%%%%%%%%%%%%%%%%%%%%%%%%%%%%%%%%%%%
%%%%% 3: Equal area and non-symmetric diffusivity %%%%%
%%%%%%%%%%%%%%%%%%%%%%%%%%%%%%%%%%%%%%%%%%%%%%%%%%%%%%%
%
\section{The equal area rule and non-symmetric diffusivity}\label{section:Nonsymmetric_diffusion}
\begin{figure}[htbp!]
    \centering
    \subcaptionbox{\label{figures:equal_area_rule_a}}{\includegraphics[width=0.45\linewidth]{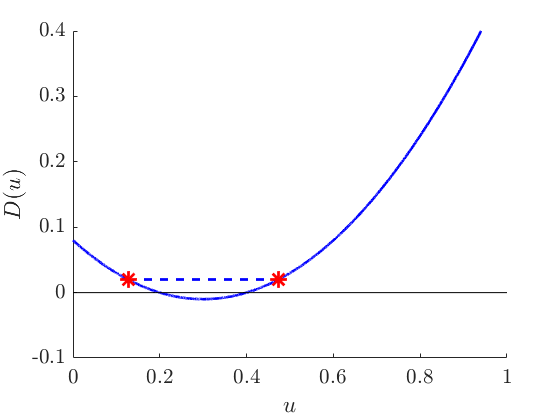}} 
    \subcaptionbox{\label{figures:equal_area_rule_b}}{\includegraphics[width=0.45\linewidth]{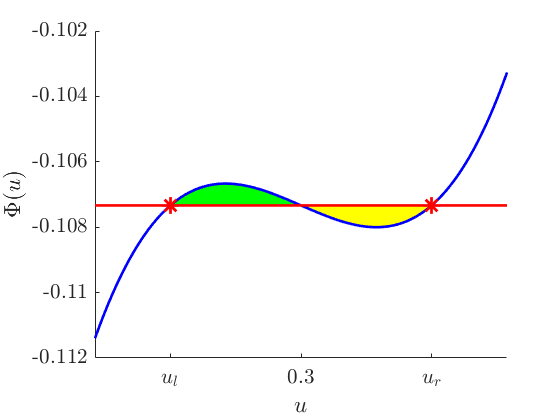}}
    \caption{(a) Stars (red) show the shock position for the solutions shown in \Cref{figures:Solution_time_dependent,figures:Travelling_wave_solution} and the dashed line (blue) shows that \(\Phi'(u) = D(u)\) is conserved across the shock. (b) The equal area rule for \(\Phi(u)\). The curve (blue) is \(\Phi(u)\) and the horizontal line (red) is the value of \(\Phi(u)\) at the shocks, \(\Phi(\ShockMin) = \Phi(\ShockMax)\). The equal area rule in \(\Phi(u)\) requires that the areas in the two lobes cut off by the horizontal line (the shaded regions) are equal.}
    \label{figures:equal_area_rule}
\end{figure}

As described in~\Cref{section:Shock}, our shock solutions are constructed by connecting two branches of an analytic multi-valued solution; here we choose the shock location according to the conditions described in equation \eqref{equation:Equal_diffusion_condition} and as shown in \Cref{figures:equal_area_rule}. In practice, the shock could be located anywhere in the multi-valued region with two conditions required to locate the shock. The first condition in \eqref{equation:Equal_diffusion_condition} requires the diffusive flux to be continuous across the shock and is a constraint used almost universally (see for example \cite{Shock-fronted:li2021,Shocks_in_nonlinear_diffusion:witelski1995}). \brev In previous works~\citep{Front_migration:pego1989,Shocks_in_nonlinear_diffusion:witelski1995, The_structure_of_internal:witelski1996}, \endrev the second condition has been chosen to correspond to the equal area in $\Phi(u)$ rule,
\begin{equation}
    \int_{\ShockMin}^{\ShockMax} \Phi(u) - \Phi(\ShockMin) \df{u}= 0 
    \label{equation:Equal_area_condition}
\end{equation}          
as illustrated in \Cref{figures:equal_area_rule_b}. The horizontal line indicates the value of $\Phi(u)$ at the shock chosen so that the two shaded regions have equal area. \citet{Front_migration:pego1989} used the equal-area rule to describe solutions of the nonlinear Cahn--Hilliard equation for phase separation, and \citet{Shocks_in_nonlinear_diffusion:witelski1995} showed how the equal-area rule arises when constructing a shock solution to the diffusion equation using a non-local regularisation
\begin{equation}
    \label{equation:regularised_diffusion}%
    \pd{u}{t} = \pd{}{x}\left(D(u)\pd{u}{x} - \varepsilon^2\pdn{u}{x}{3}\right),
\end{equation}
where \(\varepsilon \ll 1\) is a small parameter. \citet{Shock-fronted:li2021} showed that the equal-area rule also arises as a shock solution for reaction--diffusion equations using the same type of regularisation (that is, equation \eqref{equation:regularised_diffusion} with a reaction term added). Other types of regularisations, including composite regularisations, have also been used to construct shock solutions~\cite{Shock-fronted:li2021,A_geometric_singular:bradshawhajek2023}.

In~\Cref{section:ExampleSolution} the location of the shock was determined by the requirement that the \arev derivative of the flux potential \(\Phi'(u)\) must also be continuous, that is, conditions~\eqref{equation:Equal_diffusion_condition} \endrev must be satisfied across the shock. As a consequence, the diffusivity is constant across the shock (in the case of the symmetry solution presented here, relationship \eqref{equation:Reaction} means that the reaction term is also constant across the shock, but this may not be true in general). In the case of quadratic diffusivity above, the condition on the first derivative recovers the more commonly imposed equal area rule \Cref{equation:Equal_area_condition}. Here we show that this is true for any diffusivity that has two zeros and is symmetric about the midpoint of its zeros.                
        
\begin{figure}[htbp!]
    \centering
    \subcaptionbox{\label{figures:quartic_RD_D}}{\includegraphics[width=0.45\linewidth]{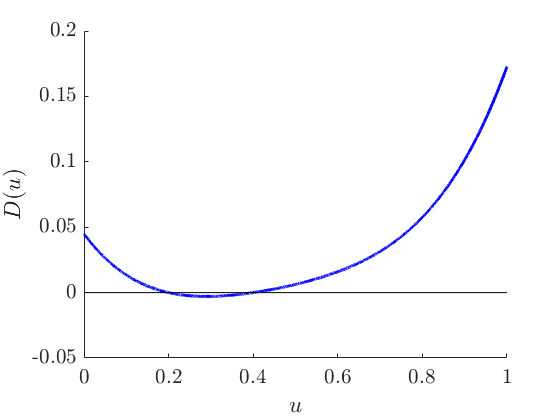}
        \begin{picture}(0,0)
            \put(-170,75){\includegraphics[width=0.23\linewidth]{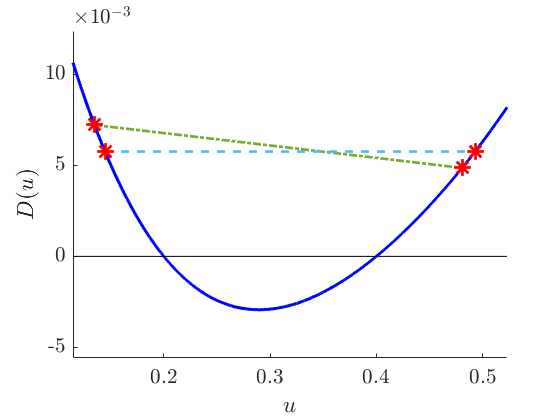}}
        \end{picture}}
    \subcaptionbox{\label{figures:quartic_RD_R}}{\includegraphics[width=0.45\linewidth]{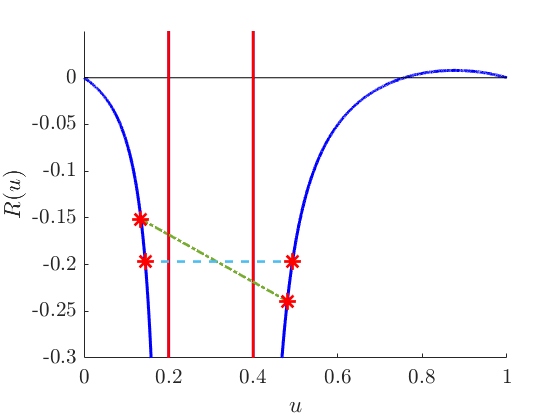}}
    \subcaptionbox{\label{figures:quartic_RD_Sol}}{\includegraphics[width=0.45\linewidth]{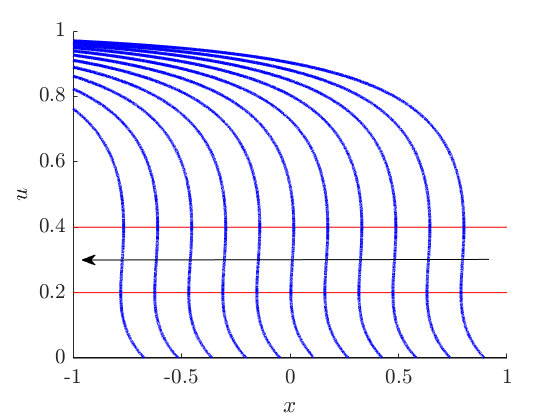}}
    \subcaptionbox{\label{figures:quartic_RD_SolZoom}}{\includegraphics[width=0.45\linewidth]{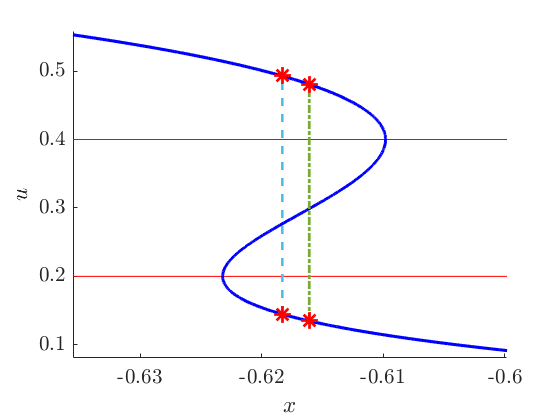}}
    \caption{Analytic solution to reaction--diffusion equation~\cref{equation:Reaction_diffusion} with non-symmetric diffusivity~\cref{equation:quarticD}. Parameter values are $\kappa = -1$, $A = 0.0448$, \(a = 0.2,\) \(b = 0.4,\) \(c = 0.6,\) and \(d = 0.2.\) (a) Diffusivity \(D(u)\)~\cref{equation:quarticD}. Inset shows the transition across the shock: horizontal dashed line (light blue) indicates the position chosen if $\Phi'(u)$ is continuous, dash-dot line (green) shows the position chosen by the equal area rule. (b) Corresponding reaction term \(R(u),\) given by~\eqref{equation:Reaction}. (c) Multi-valued travelling wave solution. Curves (blue) show the solution at different times, the arrow indicates the direction of travel. Horizontal lines (red) show where the diffusivity is zero. (d) Enlargement of the multi-valued part of the solution at a single time. Stars (red) show the end points of two possible shock positions. The dashed line (left, light blue) corresponds to \(\Phi'(u)\) is continuous. The dash-dot line (right, green) corresponds to the equal-area rule.}
    \label{figures:quartic_RD}
\end{figure}

Consider a diffusivity with two zeros such that \(D(a)=D(b)=0\) and \(D(u)<0\) for \(a<u<b\). Let \(\bar{u} = u - (a + b)/2\) so that \(D(\pm\bar{a})=0\) where \(\bar{a}=(a-b)/2\). If \(D(u)\) is symmetric about the midpoint of its zeros, \(D(\bar{u})\) is an even function. Since \(\Phi(u) = \int_{u^*}^{u}D(u') \df{u'}\), \(\bar{\Phi}=\Phi(\bar {u})\) is a vertical translation of an odd function. That is,
\begin{equation}
    f(\bar{u}) = \bar{\Phi}-\bar{\Phi}|_{\bar u=0}    
\end{equation}
is odd. Also, since \(\Phi'(\bar {u})=D(\bar {u})\) is even, the requirement that \(\Phi'(\bar{\ShockMin})=\Phi'(\bar{\ShockMax})\) is satisfied by \(\bar{\ShockMin}=-\bar{\ShockMax}\), because $\bar{\Phi}$ is a vertical translation of the odd function $f(\bar{u})$ and the requirement that \(\Phi(\bar{\ShockMin}) = \Phi(-\bar{\ShockMin})\) means that \(\Phi(\bar{\ShockMin}) = \bar{\Phi}|_{\bar u=0}\). The left-hand side of the equal area condition \eqref{equation:Equal_area_condition} then becomes
\begin{eqnarray*}
            \int_{\bar{u}_l}^{\bar{u}_r}f(\bar u)\df{\bar u} &=& \int_{\bar{u}_l}^{-\bar{u}_l}f(\bar u)\df{\bar u}\\
                &=& F(-\bar{u}_l)-F(\bar{u}_l)\\
                &=&F(\bar{u}_l)-F(\bar{u}_l)=0
\end{eqnarray*}
where \(\int f(\bar u) \df{\bar u} =F(\bar u)\), and \(F(\bar u)\) is an even function since \(f(\bar u)\) is odd. That is, the equal area rule is \arev equivalent to the requirement that \(\Phi'(u)\) be continuous, \endrev providing the diffusivity is symmetric about the midpoint of its two zeros. In contrast, if the diffusivity is not symmetric about the midpoint of its two zeros, the requirement that \arev \(\Phi'(u)\) is continuous \endrev does not recover the familiar equal area rule. For example, consider a quartic diffusivity given by
\begin{equation}
    D(u) = (u - a)(u - b)\left((u - c)^2 + d\right),
    \label{equation:quarticD}
\end{equation}
where \(a < b\), \(a,b \in (0,1)\), \(c\) arbitrary and \(d > 0\) are constants. This diffusivity has two zeros (\(D(a)=D(b)=0\)), and is not symmetric about their midpoint. \Cref{figures:quartic_RD}(a) shows this diffusivity, where the parameters have been chosen to mimic the diffusivity examined in \Cref{section:ExampleSolution}. \Cref{figures:quartic_RD}(b) shows a possible corresponding (through \Cref{equation:Reaction}) reaction-term chosen to have similar characteristics to that discussed in \Cref{section:ExampleSolution}; vertical red lines indicate the location of the zeros of the diffusivity. In both panels, the horizontal blue dashed line indicates the transition across the shock as chosen by enforcing that \arev \(\Phi'(u)\) is continuous\endrev. A direct consequence of this choice is that both \(D(u)\) and \(R(u)\) are continuous across the shock. The green dash-dot lines indicate the transition across the shock if the equal area rule is imposed. Notice that location of the shock is different, and also that neither the diffusivity nor reaction is continuous across a shock chosen in this way.
        
An implicit solution to equation \eqref{equation:Reaction_diffusion} with \eqref{equation:quarticD} and \eqref{equation:Reaction} can be constructed using the nonclassical symmetry described in \Cref{section:Introduction}. \Cref{figures:quartic_RD_Sol,figures:quartic_RD_SolZoom} show a receding travelling wave solution. The dashed line in panel (d) indicates the difference in shock location as chosen by the requirement that \arev \(\Phi'(u)\) is continuous \endrev and the equal area rule.
%
%%%%%%%%%%%%%%%%%%%%%%%%%
%%%%% 4: Discussion %%%%%
%%%%%%%%%%%%%%%%%%%%%%%%%
%
\section{Discussion and conclusion}\label{section:Discussion}
Analysis of nonlinear reaction--diffusion equations has wide-ranging application in the physical sciences. Nonclassical symmetry analysis can recover analytic solutions for specific diffusivity and reaction terms, including where diffusivity is negative. Negative diffusivity provides a phenomenological description of aggregation, and also arises in the continuum limit of discrete models with aggregation. We analyse analytic solutions to a reaction--diffusion equation with nonlinear diffusivity that is negative for some values of the density. These solutions are multi-valued, requiring the insertion of shock discontinuities. The combination of negative diffusivity and shock-fronted analytic solutions provides rich opportunities for analysis, which we explore in this work. 

\brev Exact solutions to nonlinear partial differential equations are rare. A major advantage of nonclassical symmetry techniques is their ability to uncover time-dependent exact solutions. A disadvantage of constructing solutions using a nonclassical symmetry is that we cannot impose an arbitrary initial condition. Instead, we identify the initial condition with the density profile of the symmetry solution at time \(t = 0.\) Owing to this limitation, we focus on the construction and analysis of exact symmetry solutions in this work. Furthermore, the connection between nonlinear reaction--diffusion equations with negative diffusivity and aggregation phenomena in a physical or biological system remains to be investigated. Consequently, we do not specify or solve a biologically-relevant initial-value problem for the reaction--diffusion equation with negative diffusivity in this work. Although preliminary numerical investigations using finite-difference and finite-volume methods yield solutions reminiscent of the shock-fronted solutions presented here, we do not focus on the numerical treatment of~\cref{equation:Reaction_diffusion} further. \endrev

We consider three types of analytic solutions --- a receding time-dependent solution, colliding waves, and receding travelling waves. These solutions are all sharp-fronted, such that the density \(u = 0\) at some \(x = L(t)\) with \(u'(L(t),t) \neq 0.\) In all solutions, the boundary \(L(t)\) evolves according to a Stefan-like condition. We insert shocks in the density profile to resolve regions where \(u(x,t)\) is multi-valued. We place shocks such that the flux potential \(\Phi(u)\) and its derivative \(\Phi'(u)\) are continuous across the shock. For diffusivity that is symmetric about the midpoint of its zeros (\emph{e.g.} the quadratic~\cref{equation:Diffusion}), these conditions are equivalent to the more commonly-used equal-area rule. With symmetric diffusivity, the continuity and equal-area requirements recover the shock with longest possible length. For a non-symmetric diffusivity, the largest jump coincides with the continuity condition only. \arev Since another way to write equation~\cref{equation:Reaction_diffusion} is equation~\cref{eq:Phi}, the requirement that \(\Phi'(u)\) is continuous has stronger foundation than imposing equal-area. Our planned future investigation of the entropy jump across the shock and regularisation to the PDE might yield further understanding about the link between negative diffusivity and physical phenomena. \endrev  
        
\brev The solutions presented here avoid the values of \(u(x,t)\) for which the diffusivity becomes negative. It is interesting then that the solutions do not display the usual characteristics of a standard diffusion model. This may be a consequence of the loss of mass at the moving boundary, although it is more likely that the negative values of \(D(u)\) impact the solution despite the solution not taking on the relevant values of \(u\). \endrev \arev While the travelling wave solutions presented are exact solutions that satisfy the PDE, they are not weak solutions in the usual sense of averaging through the PDE with an arbitrary compactly-supported test function. We refer the reader to~\citet{Linear_and_nonlinear_waves:whitham2011} for further information about weak solutions that contain shocks. \endrev
        
In addition to presenting the analytic solutions, we perform spectral stability analysis for constant and travelling wave solutions. We show that the constant solution \(u = 0\) is stable, but that \(u = 1\) admits long-wavelength instability for quadratic diffusivity with \(a + b > 1.\) We also prove spectral stability of travelling wave solutions to~\cref{equation:Reaction_diffusion} with quadratic diffusivity~\cref{equation:Diffusion}, for some combinations of \(a\) and \(b.\) However, open problems remain in the analysis of analytic solutions to reaction--diffusion equations with negative diffusivity. For non-symmetric diffusivity, the continuity of flux \(\Phi(u)\) and its derivative \(\Phi'(u)\) conditions for defining shocks do not yield an equal-area rule. The regularisation to the nonlinear reaction--diffusion equation~\cref{equation:Reaction_diffusion} corresponding to \arev the requirement for \(\Phi'(u)\) to be continuous \endrev remains unknown. Investigation of this and other regularisations will be the subject of future work, \arev as will a study of the jump in entropy across the shock.\endrev
\section*{Acknowledgements}
This work was supported by the Australian Research Council Discovery Program (grant numbers DP200102130 and DP230100406).
%
%%%%%%%%%%%%%%%%%%%%%%%%
%%%%% Bibliography %%%%%
%%%%%%%%%%%%%%%%%%%%%%%%
%
\bibliographystyle{abbrvnat}
\bibliography{biblo}

\begin{thebibliography}{41}
\providecommand{\natexlab}[1]{#1}
\providecommand{\url}[1]{\texttt{#1}}
\expandafter\ifx\csname urlstyle\endcsname\relax
  \providecommand{\doi}[1]{doi: #1}\else
  \providecommand{\doi}{doi: \begingroup \urlstyle{rm}\Url}\fi

\bibitem[Alt(1985)]{Models_for_mutual:alt1985}
W.~Alt.
\newblock Models for mutual attraction and aggregation of motile individuals.
\newblock In V.~Capasso, E.~Grosso, and S.~L. Paveri-Fontana, editors,
  \emph{Mathematics in Biology and Medicine}, pages 33--38, Berlin, Heidelberg,
  1985. Springer Berlin Heidelberg.

\bibitem[Aronson(1985)]{The_role_of_diffusion:aronson1985}
D.~G. Aronson.
\newblock The role of diffusion in mathematical population biology: {S}kellam
  revisited.
\newblock In V.~Capasso, E.~Grosso, and S.~L. Paveri-Fontana, editors,
  \emph{Mathematics in Biology and Medicine}, pages 2--6, Berlin, Heidelberg,
  1985. Springer Berlin Heidelberg.

\bibitem[Arrigo and Hill(1995)]{Nonclassical_symmetries:arrigo1995}
D.~J. Arrigo and J.~M. Hill.
\newblock Nonclassical symmetries for nonlinear diffusion and absorption.
\newblock \emph{Studies in Applied Mathematics}, 94\penalty0 (1):\penalty0
  21--39, 1995.

\bibitem[Bluman and Cole(1969)]{The_general_similarity:bluman1969}
G.~W. Bluman and J.~D. Cole.
\newblock The general similarity solution of the heat equation.
\newblock \emph{Journal of Mathematics and Mechanics}, 18\penalty0
  (11):\penalty0 1025--1042, 1969.

\bibitem[Bradshaw-Hajek et~al.(2023)Bradshaw-Hajek, Lizarraga, Marangell, and
  Wechselberger]{A_geometric_singular:bradshawhajek2023}
B.~H. Bradshaw-Hajek, I.~Lizarraga, R.~Marangell, and M.~Wechselberger.
\newblock A geometric singular perturbation analysis of generalised shock
  selection rules in reaction-nonlinear diffusion models.
\newblock \emph{arXiv preprint arXiv:2308.02719}, 2023.

\bibitem[Broadbridge et~al.(2002)Broadbridge, Bradshaw, Fulford, and
  Aldis]{Huxley_and_Fisher:broadbridge2002}
P.~Broadbridge, B.~H. Bradshaw, G.~R. Fulford, and G.~K. Aldis.
\newblock Huxley and {F}isher equations for gene propagation: An exact
  solution.
\newblock \emph{The Anziam Journal}, 44\penalty0 (1):\penalty0 11--20, 2002.

\bibitem[Broadbridge et~al.(2023)Broadbridge, Bradshaw-Hajek, and
  Hutchinson]{Conditionally_Integrable:broadbridge2023}
P.~Broadbridge, B.~H. Bradshaw-Hajek, and A.~J. Hutchinson.
\newblock Conditionally integrable pdes, nonclassical symmetries and
  applications.
\newblock \emph{Proceedings of the Royal Society A}, 479\penalty0
  (2276):\penalty0 20230209, 2023.

\bibitem[Courchamp et~al.(1999)Courchamp, Clutton-Brock, and
  Grenfell]{Inverse_density:courchamp1999}
F.~Courchamp, T.~Clutton-Brock, and B.~Grenfell.
\newblock Inverse density dependence and the {A}llee effect.
\newblock \emph{Trends in Ecology \& Evolution}, 14\penalty0 (10):\penalty0
  405--410, 1999.

\bibitem[Du and Lin(2010)]{Spreading-vanishing:du2010}
Y.~Du and Z.~Lin.
\newblock Spreading-vanishing dichotomy in the diffusive logistic model with a
  free boundary.
\newblock \emph{SIAM Journal on Mathematical Analysis}, 42\penalty0
  (1):\penalty0 377--405, 2010.

\bibitem[Edwards et~al.(2018)Edwards, Bradshaw-Hajek, Munoz-Lopez, Waterhouse,
  and Anderssen]{Compactly_Supported_Solutions:edwards2018}
M.~P. Edwards, B.~H. Bradshaw-Hajek, M.~J. Munoz-Lopez, P.~M. Waterhouse, and
  R.~S. Anderssen.
\newblock Compactly supported solutions of reaction--diffusion models of
  biological spread.
\newblock In \emph{Agriculture as a Metaphor for Creativity in All Human
  Endeavors}, pages 125--138. Springer, 2018.

\bibitem[El-Hachem et~al.(2019)El-Hachem, McCue, Jin, Du, and
  Simpson]{Revisiting_the_Fisher:el-hachem2019}
M.~El-Hachem, S.~W. McCue, W.~Jin, Y.~Du, and M.~J. Simpson.
\newblock Revisiting the {F}isher--{K}olmogorov--{P}etrovsky--{P}iskunov
  equation to interpret the spreading--extinction dichotomy.
\newblock \emph{Proceedings of the Royal Society A: Mathematical, Physical and
  Engineering Sciences}, 475\penalty0 (2229):\penalty0 20190378, 2019.

\bibitem[El-Hachem et~al.(2021)El-Hachem, McCue, and
  Simpson]{Invading_and_receding:el-hachem2021}
M.~El-Hachem, S.~W. McCue, and M.~J. Simpson.
\newblock Invading and receding sharp-fronted travelling waves.
\newblock \emph{Bulletin of Mathematical Biology}, 83\penalty0 (4):\penalty0
  1--25, 2021.

\bibitem[Fadai(2021)]{Semi-Infinite:fadai2021}
N.~T. Fadai.
\newblock Semi-infinite travelling waves arising in a general
  reaction–diffusion {{Stefan}} model.
\newblock \emph{Nonlinearity}, 34\penalty0 (2):\penalty0 725--743, 2021.
\newblock ISSN 0951-7715.
\newblock \doi{10.1088/1361-6544/abd07b}.

\bibitem[Ferracuti et~al.(2009)Ferracuti, Marcelli, and
  Papalini]{Travelling_waves_in_some:ferracuti2009}
L.~Ferracuti, C.~Marcelli, and F.~Papalini.
\newblock Travelling waves in some reaction-diffusion-aggregation models.
\newblock \emph{Advances in Dynamical Systems and Applications}, 4\penalty0
  (1):\penalty0 19--33, 2009.

\bibitem[Goard(2008)]{Goard2008}
J.~Goard.
\newblock Finding symmetries by incorporating initial conditions as side
  conditions.
\newblock \emph{European Journal of Applied Mathematics}, 19\penalty0
  (6):\penalty0 701--715, Dec. 2008.
\newblock ISSN 0956-7925, 1469-4425.
\newblock \doi{10.1017/S0956792508007705}.

\bibitem[Goard and Broadbridge(1996)]{Nonclassical_symmetry:goard1996}
J.~Goard and P.~Broadbridge.
\newblock Nonclassical symmetry analysis of nonlinear reaction-diffusion
  equations in two spatial dimensions.
\newblock \emph{Nonlinear Analysis: Theory, Methods \& Applications},
  26\penalty0 (4):\penalty0 735--754, 1996.

\bibitem[Gr{\"u}nbaum and Okubo(1994)]{Modelling_social:grunbaum1994}
D.~Gr{\"u}nbaum and A.~Okubo.
\newblock Modelling social animal aggregations.
\newblock In S.~A. Levin, editor, \emph{Frontiers in Mathematical Biology},
  pages 296--325, Berlin, Heidelberg, 1994. Springer Berlin Heidelberg.

\bibitem[Harley et~al.(2014)Harley, van Heijster, Marangell, Pettet, and
  Wechselberger]{Existence_of_traveling_wave:harley2014}
K.~Harley, P.~van Heijster, R.~Marangell, G.~J. Pettet, and M.~Wechselberger.
\newblock Existence of traveling wave solutions for a model of tumor invasion.
\newblock \emph{SIAM Journal on Applied Dynamical Systems}, 13\penalty0
  (1):\penalty0 366--396, 2014.

\bibitem[Ibragimov(1999)]{Elementary_Lie:ibragimov1999}
N.~H. Ibragimov.
\newblock \emph{Elementary Lie group analysis and ordinary differential
  equations}, volume 197.
\newblock Wiley New York, 1999.

\bibitem[Johnston et~al.(2017)Johnston, Baker, McElwain, and
  Simpson]{Co-operation_competition:johnston2017}
S.~T. Johnston, R.~E. Baker, D.~L.~S. McElwain, and M.~J. Simpson.
\newblock Co-operation, competition and crowding: a discrete framework linking
  {A}llee kinetics, nonlinear diffusion, shocks and sharp-fronted travelling
  waves.
\newblock \emph{Scientific Reports}, 7\penalty0 (1):\penalty0 1--19, 2017.

\bibitem[Kapitula and Promislow(2013)]{Spectral_and_dynamical:kapitula2013}
T.~Kapitula and K.~Promislow.
\newblock \emph{Spectral and dynamical stability of nonlinear waves}, volume
  457.
\newblock Springer, 2013.

\bibitem[Keller and Segel(1970)]{Keller1970}
E.~F. Keller and L.~A. Segel.
\newblock Initiation of slime mold aggregation viewed as an instability.
\newblock \emph{J. Theor. Biol.}, 26\penalty0 (3):\penalty0 399--415, 1970.
\newblock ISSN 0022-5193.
\newblock \doi{10.1016/0022-5193(70)90092-5}.

\bibitem[Kuzmin and Ruggerini(2011)]{Front_propagation:kuzmin2011}
M.~Kuzmin and S.~Ruggerini.
\newblock Front propagation in diffusion-aggregation models with bi-stable
  reaction.
\newblock \emph{Discrete \& Continuous Dynamical Systems --- B}, 16\penalty0
  (3):\penalty0 819--833, 2011.

\bibitem[Lewis and Kareiva(1993)]{Allee_dynamics:lewis1993}
M.~A. Lewis and P.~Kareiva.
\newblock Allee dynamics and the spread of invading organisms.
\newblock \emph{Theoretical Population Biology}, 43\penalty0 (2):\penalty0
  141--158, 1993.

\bibitem[Li et~al.(2020)Li, van Heijster, Marangell, and
  Simpson]{Travelling_wave_solutions:li2020}
Y.~Li, P.~van Heijster, R.~Marangell, and M.~J. Simpson.
\newblock Travelling wave solutions in a negative nonlinear diffusion--reaction
  model.
\newblock \emph{Journal of Mathematical Biology}, 81\penalty0 (6-7):\penalty0
  1495--1522, 2020.

\bibitem[Li et~al.(2021)Li, van Heijster, Simpson, and
  Wechselberger]{Shock-fronted:li2021}
Y.~Li, P.~van Heijster, M.~J. Simpson, and M.~Wechselberger.
\newblock Shock-fronted travelling waves in a reaction--diffusion model with
  nonlinear forward--backward--forward diffusion.
\newblock \emph{Physica D: Nonlinear Phenomena}, 423:\penalty0 132916, 2021.

\bibitem[Lie(1880)]{Lie_article}
S.~Lie.
\newblock Theorie der transformationsgruppen {I}.
\newblock \emph{Mathematische Annalen}, 16\penalty0 (4):\penalty0 441--528,
  1880.

\bibitem[Lundberg and Totik(2013)]{Lundberg2013}
E.~Lundberg and V.~Totik.
\newblock Lemniscate growth.
\newblock \emph{Anal. Math. Phys.}, 3\penalty0 (1):\penalty0 45--62, 2013.
\newblock \doi{10.1007/s13324-012-0038-1}.

\bibitem[Maini et~al.(2006)Maini, Malaguti, Marcelli, and
  Matucci]{Diffusion_aggregation:maini2006}
P.~K. Maini, L.~Malaguti, C.~Marcelli, and S.~Matucci.
\newblock Diffusion-aggregation processes with mono-stable reaction terms.
\newblock \emph{Discrete \& Continuous Dynamical Systems --- B}, 6\penalty0
  (5):\penalty0 1175--1189, 2006.

\bibitem[Moitsheki et~al.(2005)Moitsheki, Broadbridge, and
  Edwards]{Symmetry_solutions:moitsheki2005}
R.~J. Moitsheki, P.~Broadbridge, and M.~P. Edwards.
\newblock Symmetry solutions for transient solute transport in unsaturated
  soils with realistic water profile.
\newblock \emph{Transport in Porous Media}, 61:\penalty0 109--125, 2005.

\bibitem[Olver(2000)]{Applications_of_Lie:olver1993}
P.~J. Olver.
\newblock \emph{Applications of Lie groups to differential equations}, volume
  107.
\newblock Springer Science \& Business Media, 2000.

\bibitem[Pego(1989)]{Front_migration:pego1989}
R.~L. Pego.
\newblock Front migration in the nonlinear {C}ahn-{H}illiard equation.
\newblock \emph{Proceedings of the Royal Society of London. A. Mathematical and
  Physical Sciences}, 422\penalty0 (1863):\penalty0 261--278, 1989.

\bibitem[Pettet et~al.(2000)Pettet, McElwain, and
  Norbury]{Lokta-Volterra_equations:pettet2000}
G.~J. Pettet, D.~L.~S. McElwain, and J.~Norbury.
\newblock Lotka-{V}olterra equations with chemotaxis: walls, barriers and
  travelling waves.
\newblock \emph{Mathematical Medicine and Biology: A Journal of the IMA},
  17\penalty0 (4):\penalty0 395--413, 2000.

\bibitem[Popescu and Dietrich(2004)]{Model_for_spreading:popescu2004}
M.~N. Popescu and S.~Dietrich.
\newblock Model for spreading of liquid monolayers.
\newblock \emph{Physical Review E}, 69\penalty0 (6):\penalty0 061602, 2004.

\bibitem[Sherratt and Murray(1990)]{Models_of_epidermal:sherratt1990}
J.~A. Sherratt and J.~D. Murray.
\newblock Models of epidermal wound healing.
\newblock \emph{Proceedings of the Royal Society of London. Series B:
  Biological Sciences}, 241\penalty0 (1300):\penalty0 29--36, 1990.

\bibitem[Stephens and Sutherland(1999)]{Consequences_of_the_Allee:stephens1999}
P.~A. Stephens and W.~J. Sutherland.
\newblock Consequences of the {A}llee effect for behaviour, ecology and
  conservation.
\newblock \emph{Trends in ecology \& evolution}, 14\penalty0 (10):\penalty0
  401--405, 1999.

\bibitem[Turchin(1989)]{Population_consequences:turchin1989}
P.~Turchin.
\newblock Population consequences of aggregative movement.
\newblock \emph{Journal of Animal Ecology}, 58\penalty0 (1):\penalty0 75--100,
  1989.

\bibitem[Wechselberger and Pettet(2010)]{Folds_canards:wechselberger2010}
M.~Wechselberger and G.~J. Pettet.
\newblock Folds, canards and shocks in advection--reaction--diffusion models.
\newblock \emph{Nonlinearity}, 23\penalty0 (8):\penalty0 1949, 2010.

\bibitem[Whitham(2011)]{Linear_and_nonlinear_waves:whitham2011}
G.~B. Whitham.
\newblock \emph{Linear and nonlinear waves}.
\newblock John Wiley \& Sons, 2011.

\bibitem[Witelski(1995)]{Shocks_in_nonlinear_diffusion:witelski1995}
T.~P. Witelski.
\newblock Shocks in nonlinear diffusion.
\newblock \emph{Applied Mathematics Letters}, 8\penalty0 (5):\penalty0 27--32,
  1995.

\bibitem[Witelski(1996)]{The_structure_of_internal:witelski1996}
T.~P. Witelski.
\newblock The structure of internal layers for unstable nonlinear diffusion
  equations.
\newblock \emph{Studies in Applied Mathematics}, 97\penalty0 (3):\penalty0
  277--300, 1996.

\end{thebibliography}
\end{document}